\documentclass[journal]{IEEEtran}
\pdfoutput=1
\usepackage{subfigure}
\usepackage{subfigmat}
\usepackage{amsmath}
 \usepackage{amssymb}
 \usepackage{graphicx}                                                      
\ifCLASSINFOpdf
\else
\fi
\begin{document}

\title{Cooperative Robot Control and Concurrent Synchronization of Lagrangian Systems}
%
%
%

\author{Soon-Jo Chung,~\IEEEmembership{Member,~IEEE,} and Jean-Jacques Slotine
\thanks{Assistant Professor of Aerospace Engineering, Iowa State University, Ames, IA 50011, E-mail: sjchung@alum.mit.edu.}
\thanks{Professor of Mechanical Engineering \& Information Sciences, Professor of Brain \& Cognitive
Sciences, Massachusetts Institute of Technology, Cambridge, MA 02319, E-mail: jjs@mit.edu.}}


%
%

\markboth{}{Chung, Slotine Cooperative Robot Control and Concurrent Synchronization of Lagrangian Systems}
%



\maketitle

\newtheorem{theorem}{\textbf{Theorem}}
\newtheorem{lemma}{\textbf{Lemma}}
\newtheorem{proposition}[theorem]{Proposition}
\newtheorem{corollary}{\textbf{Corollary}}
\newtheorem{remark}{Remark}

\newenvironment{definition}[1][Definition]{\begin{trivlist}
\item[\hskip \labelsep {\bfseries #1}]}{\end{trivlist}}
\newenvironment{example}[1][Example]{\begin{trivlist}
\item[\hskip \labelsep {\bfseries #1}]}{\end{trivlist}}

\newcommand{\qed}{\nobreak \ifvmode \relax \else
      \ifdim\lastskip<1.5em \hskip-\lastskip
      \hskip1.5em plus0em minus0.5em \fi \nobreak
      \vrule height0.75em width0.5em depth0.25em\fi}

%
\begin{abstract}
Concurrent synchronization is a regime where diverse groups of fully synchronized dynamic systems stably coexist. We study global exponential synchronization and concurrent synchronization in the context of Lagrangian systems control.  In a network constructed by adding diffusive couplings to robot manipulators or mobile robots, a decentralized tracking control law globally exponentially synchronizes an arbitrary number of robots, and represents a generalization of the average consensus problem. Exact nonlinear stability guarantees and synchronization conditions are derived by contraction analysis. The proposed decentralized strategy is further extended to adaptive synchronization and partial-state coupling.
\end{abstract}
\section{Introduction}
Distributed and decentralized synchronization phenomena of large
dynamic groups are areas of intense research. In this article,
synchronization is defined as a complete match of all configuration
variables of each dynamic system such that
$\mathbf{x}_1=\mathbf{x}_2=\cdots=\mathbf{x}_p$ and $p$ denotes the
number of sub-systems in the network. While such a definition directly
concerns the angular position synchronization problem, this paper also
addresses the synchronization of biased variables with application to
the translational position coordination problem. In the latter case,
synchronization corresponds to
$\mathbf{y}_1=\mathbf{y}_2=\cdots=\mathbf{y}_p$ where $\mathbf{y}_i$,
$1\leq i\leq p$ connotes a vector of biased variables constructed from
the configuration vector $\mathbf{x}_i$ such that
$\mathbf{y}_i(t)=\mathbf{x}_i(t)+\mathbf{b}_i(t)$ and the separation
vector $\mathbf{b}_i(t)$ is independent of the dynamics. Furthermore, we
construct complex robot networks where multiple groups of fully
synchronized elements coexist. Such \emph{concurrent synchronization}
seems pervasive in biology, and in particular in the brain
where multiple rhythms coexist and neurons can exhibit many
qualitatively different types of
oscillations~\cite{Ref:contraction_sync}.

The objective of this paper is to establish a unified
synchronization framework that can achieve both the synchronization of the configuration variables of the robots and the stable tracking of a common desired trajectory. Although an uncoupled trajectory tracking control law, in the absence of external disturbances, would achieve synchronization to a common desired trajectory, non-identical disturbances justify the mutual synchronization of the system variables. On the other hand, the synchronization to the average of initial conditions is not sufficient for multi-robot or multi-vehicle systems where a desired trajectory is explicitly given. For example, a large
swarm of robots can first synchronize their attitudes and positions to form a certain formation
pattern, then track the common desired trajectory to accomplish the given
mission. In production processes, such as manufacturing and automotive applications, where high flexibility, manipulability,
and maneuverability cannot be achieved by a single system~\cite{Ref:robot_sync}, there has been widespread interest in cooperative schemes for multiple robot manipulators that track a predefined trajectory. A stellar formation
flight interferometer~\cite{Ref:tether,Ref:FormationFlying_CHUNG} is another example where
precision control of relative spacecraft motions is indispensable. The proposed synchronization tracking control law can be implemented
for such purposes, where a common desired trajectory can be explicitly given. The proposed strategy can achieve more efficient and robust
performance through local interactions, especially in the presence of non-identical external disturbances. Further, we generalize the proposed control law such that multiple dynamic systems can synchronize themselves from arbitrary initial conditions without the need for a common reference trajectory. As a result, other potential applications
include oscillation synchronization of robotic
locomotion~\cite{Ref:sync_Pitti,SeoA}, and
tele-manipulation of
robots~\cite{Ref:Spong_network,Ref:tele}.

The main contributions of this work can be stated as follows.
\begin{itemize}
\item Concurrent synchronization that exploits the multiple time scale behaviors from two types of inputs (a reference trajectory and local couplings) permits construction of a complex time-varying network comprised of numerous heterogeneous systems.
\item In contrast with prior work on consensus and flocking problems using graphs,
the proposed strategy primarily deals with dynamic networks
consisting of nonlinear time-varying dynamics.
 \item We use contraction
analysis~\cite{Ref:contraction1,Ref:contraction3} as our main nonlinear stability tool, thereby
deriving exact and global results with exponential convergence, as opposed to asymptotic convergence of prior work.
  \item The
proposed control laws are of a decentralized form requiring only
local velocity/position coupling feedback for global exponential
convergence, thereby facilitating implementation in real systems.
  \item The theory is generalized and extended to multi-robot
systems with non-identical dynamics, linear coupling control, partial state coupling,
uni-directional coupling, and adaptive control.
\end{itemize}
\subsection{Comparison with Related Work}\label{sec:sync_previouswork}
The consensus problems on graph~\cite{Ref:Murray} and the
coordination of multi-agent
systems~\cite{Jadbabaie,Ref:mobile,CLF:2002,Ref:multi_coord} are
closely related with the synchronization problem. In particular, the
use of graph theory and Laplacian produced many interesting
results~\cite{Jadbabaie,Ref:passive_decomp2,Ref:mobile,Ref:Mesbahi,Ref:Murray,Ref:WeiRen1,Ref:WeiRen2}.
However, the synchronization to the average of initial conditions might not be directly applicable to multi-robot and multi-vehicle systems, where a desired trajectory is explicitly defined. A recent work~\cite{Ref:WeiRen1} studied the consensus problem with a time-varying reference state, based on a single integrator model. In essence, the aforementioned work mainly
deals with simple dynamic models such as linear systems and
single or double integrator models without nonlinearly coupled inertia matrices. In contrast, we aim at addressing
highly nonlinear systems (e.g. helicopters, attitude
dynamics of spacecraft, walking robots, and manipulator robots). As
shall be seen later, the proof of the synchronization for network
systems that possess nonlinearly coupled inertia matrices is more
involved. This paper focuses on such dynamic networks
consisting of highly nonlinear systems.

One notable work~\cite{Ref:robot_sync}
introduced a nonlinear tracking control law that synchronized multiple
robot manipulators in order to track a common desired trajectory. Due to its all-to-all coupling requirement, the number of variables to be estimated increases with
the number of robots to be synchronized, which imposes a significant
communication burden. Additionally, the feedback of estimated
acceleration errors requires unnecessary information and complexity.
Thus, a method to eliminate both the all-to-all coupling and the
feedback of the acceleration terms is explored in this paper. Another nice approach to the synchronization of robot networks is to
exploit the passivity of the input-output
dynamics~\cite{Ref:Spong_network,Ref:Spong_Sync}. Its property of
robustness to time delays is particularly attractive, while robot
dynamics are passive only with velocity outputs unless composite
variables are employed. In addition, the mutual synchronization
problem, which not only synchronizes the sub-members but also
enforces them to follow a common reference trajectory, is not
addressed. In particular, it shall be shown that our proposed control law generalizes the robot control law presented in \cite{Ref:Spong_Sync}, while our convergence result is stronger (exponential). One recent work~\cite{IhleArcak} used the passivity property for the path-following system to synchronize the path variables. The passive decomposition~\cite{Ref:passive_decomp1,Ref:passive_decomp2} describes a
strategy of decoupling into the internal group formation shape and the total group maneuver. Due to its dependency on a centralized
control architecture, the decoupling is not generally ensured under
the decentralized control. \cite{Ref:passive_decomp2} considers
linear double integrator models. Another work~\cite{Ref:robot_formation} proposed a nice control framework for controlling and coordinating a group of nonholonomic mobile robots, although the exact synchronization of individual vehicles was not pursued. One recent work~\cite{PaleyLeonard} presented a framework called particles with coupled oscillator dynamics so that collective spatial patterns emerge. Our recent work~\cite{Ref:FormationFlying_CHUNG} also discussed more complex spacecraft dynamics to achieve the phase synchronization of spiral or circular translational trajectories in three dimensions as well as the synchronization of rotational dynamics. Also, it is important to note that the adaptive synchronization of multiple robotic manipulators in a single local coupling configuration was studied in~\cite{Sun1,Sun2}.

The proposed synchronization framework
using contraction analysis, first reported in \cite{Ref:sjchung_CDC}, has a clear advantage in its broad
applications to a larger class of identical or nonidentical
nonlinear systems even with complex coupling geometry including uni-directional couplings and
partial degrees-of-freedom couplings, non-passive input-output,
and time delays~\cite{Ref:sjchung_PhD,Ref:contraction_timedelay}, while ensuring a simple
decentralized coupling control law (see Figure \ref{robot_diagram} for
network structures permitted here). The proofs and new results from \cite{Ref:sjchung_CDC} are expanded in this paper. In particular, concurrent synchronization that exploits two different types of inputs has not been studied in the literature.
\subsection{Organization}
Section
\ref{Sec:modelingCH5} describes modeling of robots
based on the Lagrangian formulation, and summarizes the key
stability theorems. While Section~\ref{sec:sync} presents the main control law and its tracking stability, the proof of exponential synchronization is more involved
and treated separately in
Section \ref{sec:synsync}. The remainder of the paper further
highlights the unique contributions of this work. Section \ref{sec:concurrent} elucidates the concurrent
synchronization of complex networks and the leader-follower problem. The main idea of this paper is extended to
linear Proportional-Derivative (PD) couplings, and limited partial-state couplings in Section
\ref{sec:ext_sim}. Key simulation results are also presented in Section
\ref{sec:ext_sim}.
\section{Modeling and Nonlinear Stability Tools for Multi-Robot Networks}\label{Sec:modelingCH5}
\subsection{Lagrangian Systems}\label{Sec:Lagrangian}
This paper is devoted
to the use of the Lagrangian formulation for its simplicity in
dealing with complex systems involving multiple dynamics. The
equations of motion for a robot with multiple joints ($\mathbf{q}_i
\in \mathbb{R}^n$) can be derived by exploiting the Euler-Lagrange
equations:
\begin{equation}\label{lagrangeCH5}
L_i=\frac{1}{2}{\mathbf{\dot{q}}_i}^{T}\mathbf{M}_i(\mathbf{q}_i)\mathbf{\dot{q}}_i-V_i, \ \ \frac{d}{dt}\frac{\partial
L_i(\mathbf{q}_i,\mathbf{\dot{q}}_i)}{\partial {\mathbf{\dot
q}}_i}-\frac{\partial L_i(\mathbf{q}_i,\mathbf{\dot{q}}_i)}{\partial
\mathbf{q}_i}=\mathbf{\tau}_i
\end{equation}
where $i$, $(1\leq i \leq p)$ denotes the index of robots or
dynamic systems comprising a network, and $p$ is the total number
of the individual elements. Equation (\ref{lagrangeCH5}) can be
represented as
\begin{align}\label{NL_single_compactCH5}
  \mathbf{M}_i(\mathbf{q}_i)\mathbf{\ddot{q}}_i&+
       \mathbf{C}_i(\mathbf{q}_i,\mathbf{\dot{q}}_i)\mathbf{\dot{q}}_i+\mathbf{g}_i(\mathbf{q}_i)
      =\mathbf{\tau}_i
\end{align}
where $\mathbf{g}_i(\mathbf{q}_i)=\frac{\partial V_i}{\mathbf{q}_i}$,
and, $\mathbf{\tau}_i$ is a generalized force or torque acting on
the $i$-th robot.

Note that we define $\mathbf{C}_i(\mathbf{q}_i,\mathbf{\dot{q}}_i)$
such that $(\mathbf{\dot{M}}_i-2\mathbf{C}_i)$ is
skew-symmetric~\cite{Ref:Slotine}, and this property plays a central
role in our stability analysis using contraction
theory~\cite{Ref:sjchung_PhD}.

The following key assumptions are used throughout this paper. Since the main applications of the present paper include fully actuated robot manipulators and spacecraft~\cite{Ref:FormationFlying_CHUNG}, the robot system in (\ref{NL_single_compactCH5}) is fully  actuated. In
other words, the number of control inputs is equal to the dimension
of their configuration manifold ($=n$). The mass-inertia matrix
$\mathbf{M}(\mathbf{q})$ is assumed to be uniformly positive definite, for all
positions $\mathbf{q}$ in the robot workspace~\cite{Ref:Slotine}.
\subsection{Contraction Analysis for Global and Exponential Stability}
\noindent\emph{(1) Modular Stability Analysis:} Although one popular method for modular stability analysis is to exploit the passivity formalism~\cite{Ref:Slotine}, we use contraction theory~\cite{Ref:contraction1,Ref:contraction_sync,Ref:contraction2,Ref:contraction3} as an alternative tool for analyzing modular stability of coupled nonlinear systems. In particular, contraction analysis has more general and intuitive combination properties (e.g., hierarchies) than the passivity method, since it involves a state-space rather than an input-output method.\\
\emph{(2) Differential State-State Analysis:} Lyapunov's linearization method
indicates that the local stability of the nonlinear system can be
analyzed using its differential approximation. What is new in
contraction theory is that a differential stability analysis can be
made exact, thereby yielding global results on the time-varying nonlinear system.\\
\emph{(3) Stronger Stability Results (Exponential and Global):} For a robot dynamic model in (\ref{NL_single_compactCH5}) and a time-varying tracking control law, the straightforward use of Barbalat's lemma or LaSalle-Yoshizawa theorem~\cite{Krstic} yields asymptotic convergence results. While global exponential stability can be proven by using a Lyapunov function with a cross-term and additional constraints~\cite{Ref:Arimoto1,Ref:Arimoto2,Khalil:2002}, such a method is ad hoc, as compared to contraction theory. While exact exponential convergence might not be achievable in real systems due to the modeling errors, we believe that finding an explicit convergence rate with exponential stability, if possible, is important due to its superior tracking performance and property of robustness with respect to perturbations (e.g., see p. 339--350 in~\cite{Khalil:2002}).

A brief review of the results from
\cite{Ref:contraction1,Ref:contraction2,Ref:contraction3} is
presented in this section. Note
that contraction theory is a generalization of the classical
Krasovskii's theorem~\cite{Ref:Slotine}, and that approaches closely
related to contraction, although not based on differential analysis,
can be traced back to \cite{Demi,Hart} and even to \cite{Lew}.

Consider a smooth nonlinear system
\begin{equation}\label{xfx}
{\mathbf{\dot
x}}(t)=\mathbf{f}(\mathbf{x}(t),\mathbf{u}(\mathbf{x},t),t)
\end{equation}
where $\mathbf{x}(t)\in\mathbb{R}^n$, and $\mathbf{f}:
\mathbb{R}^n\times\mathbb{R}^m\times\mathbb{R}_{+}\rightarrow\mathbb{R}^n$.
A virtual displacement, $\delta\mathbf{x}$ is defined as an
infinitesimal displacement at a fixed time-- a common supposition in
the calculus of variations.
\begin{theorem}\label{Thm:contraction}
\emph{For the system in (\ref{xfx}), if there exists a uniformly positive
definite metric,
\begin{equation}\label{metric}
\mathbf{M}(\mathbf{x},t)={\mathbf{\Theta}}(\mathbf{x},t)^{T}{\mathbf{\Theta}}(\mathbf{x},t)
\end{equation}
where $\mathbf{\Theta}$ is some smooth coordinate transformation of
the virtual displacement,
$\delta\mathbf{z}={\mathbf{\Theta}}\delta\mathbf{x}$, such that the
associated generalized Jacobian, $\bf F$ is uniformly negative
definite, i.e., $\exists \lambda >0$ such that
\begin{equation}\label{jacobian}
\mathbf{F}=\left(\mathbf{\dot{\Theta}}{(\mathbf{x},t)}+{\mathbf{\Theta}(\mathbf{x},t)}\frac{\partial
\mathbf{f}}{\partial
\mathbf{x}}\right){\mathbf{\Theta}(\mathbf{x},t)}^{-1} \le - \lambda
{\bf I},
\end{equation}
then all system trajectories converge globally to a single
trajectory exponentially fast regardless of the initial conditions,
with a global exponential convergence rate of the largest
eigenvalues of the symmetric part of $\mathbf{F}$.}
\end{theorem}
Such a system is said to be contracting. The proof is given in
\cite{Ref:contraction1}. Equivalently, the system is contracting if
$\exists \lambda >0$ such that
\begin{equation}\label{length}
\mathbf{\dot{M}}+\left(\frac{\partial \mathbf{f}}{\partial
\mathbf{x}}\right)^T\mathbf{{M}}+\mathbf{{M}}\frac{\partial
\mathbf{f}}{\partial \mathbf{x}} \le - 2 \lambda {\bf M}
\end{equation}
Equation (\ref{length}) is useful for the stability proof of a Lagrangian system, since the inertia matrix $\mathbf{M}(\mathbf{q})$ of the robot dynamics in (\ref{NL_single_compactCH5}) can be chosen as the metric $\mathbf{{M}}$ in (\ref{length}).

The key advantage of contraction analysis is its superior combination property as follows.
\subsection{Contraction of Coupled Systems}
The following theorems are used to derive stability and
synchronization of coupled Lagrangian systems.
\begin{theorem}\label{Thm:hierc}
\emph{Hierarchical
combination~\cite{Ref:contraction2}. Consider two
contracting systems, of possibly different dimensions and metrics,
and connect them in series, leading to a smooth virtual dynamics of
the form
\begin{equation*}
\frac{d}{dt} \begin{pmatrix} \delta\mathbf{z_1}\\
\delta\mathbf{z_2}
\end{pmatrix}=
\begin{pmatrix} \mathbf{F_{11}} & \mathbf{0}
\\ \mathbf{F_{21}} & \mathbf{F_{22}} \end{pmatrix} \begin{pmatrix} \delta\mathbf{z_1}\\
\delta\mathbf{z_2}
\end{pmatrix}
\end{equation*}
Then, the combined system is contracting if $\mathbf{F_{21}}$ is
bounded.}
\end{theorem}
\begin{proof} see~\cite{Ref:sjchung_PhD,Ref:contraction2,Ref:contraction5}.\end{proof}
\begin{theorem}\label{Thm:partial}
\emph{Partial contraction~\cite{Ref:contraction3}. Consider a nonlinear system of the form ${\mathbf{ \dot
x}}=\mathbf{f}(\mathbf{x},\mathbf{x},t)$ and assume that the
auxiliary system ${\mathbf{\dot
y}}=\mathbf{f}(\mathbf{y},\mathbf{x},t)$ is contracting with respect
to $\mathbf{y}$. If a particular solution of the auxiliary
$\mathbf{y}$-system verifies a specific smooth property, then all
trajectories of the original x-system verify this property
exponentially. The original system is said to be partially
contracting.}
\end{theorem}
\begin{proof}See \cite{Ref:contraction3} for the virtual observer-like $\mathbf{y}$ system.\end{proof}
\begin{theorem}\label{Thm:sync}\emph{Synchronization~\cite{Ref:contraction3}. Consider two coupled systems. If the dynamics equations verify
\begin{equation*}
{\mathbf{\dot x}}_1-\mathbf{f}(\mathbf{x}_1,t)={\mathbf{\dot
x}}_2-\mathbf{f}(\mathbf{x}_2,t)
\end{equation*}
where the function $\mathbf{f}(\mathbf{x},t)$ is contracting in an
input-independent metric, then $\mathbf{x}_1$ and $\mathbf{x}_2$
will converge to each other exponentially, regardless of the initial
conditions. Mathematically, stable concurrent synchronization
corresponds to convergence to a flow-invariant linear subspace of
the global state space~\cite{Ref:contraction_sync}.}
\end{theorem}
\begin{proof}
This can be proven by constructing the virtual system ${\mathbf{\dot y}}-\mathbf{f}(\mathbf{y},t)=\mathbf{u}(t)$ and using Theorem~\ref{Thm:partial} (see \cite{Ref:contraction3}).
\end{proof}
\begin{remark}
Whereas Theorem~\ref{Thm:hierc} can be proven for different metrics, Theorem~\ref{Thm:sync} requires the same metric (e.g., inertia matrix) among the coupled systems. Hence, Theorem~\ref{Thm:sync} cannot be directly applied to coupled Lagrangian systems. This is one of the motivations of the current paper, and elucidated in Section~\ref{sec:synsync}.\end{remark}
\section{Control Law and Its Tracking Stability}\label{sec:sync}
A tracking controller introduced in
this section achieves not only global exponential
synchronization of the configuration variables, but also global
exponential convergence to the desired trajectory.
\nopagebreak
\subsection{Proposed Synchronization Control Strategy}\label{sec:tack_sync}
\begin{figure}[t]
 \begin{subfigmatrix}{2}
 \subfigure[with bi-directional couplings]{\includegraphics{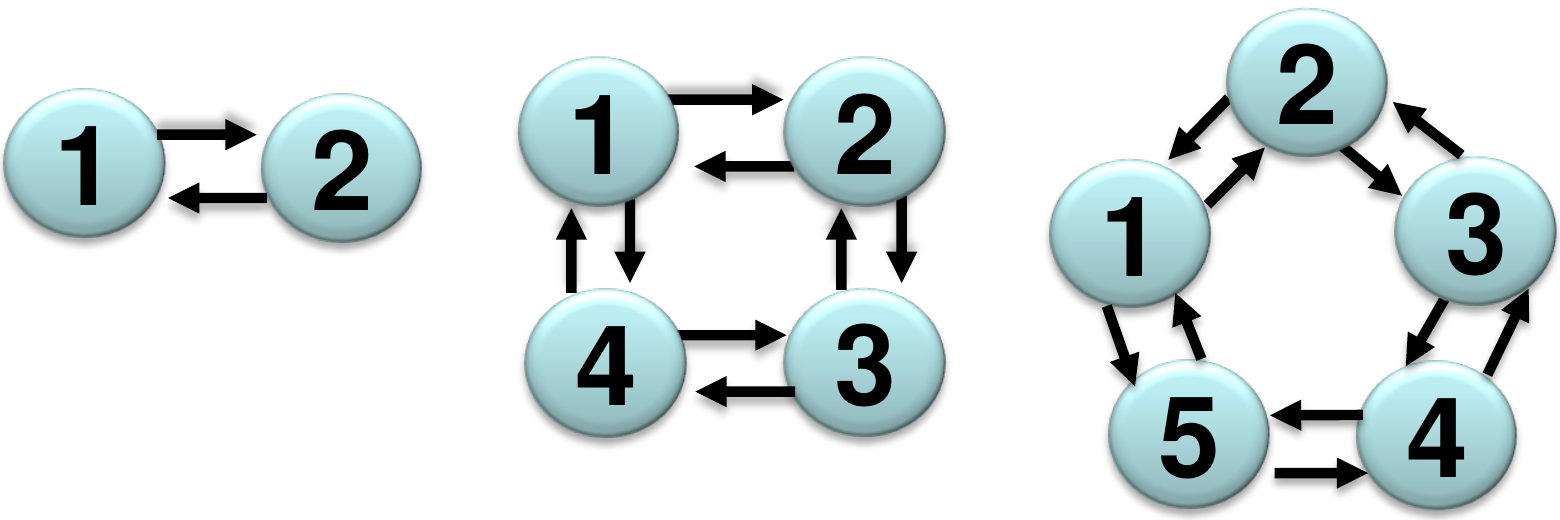}}
 \subfigure[heterogeneous robots]{\includegraphics{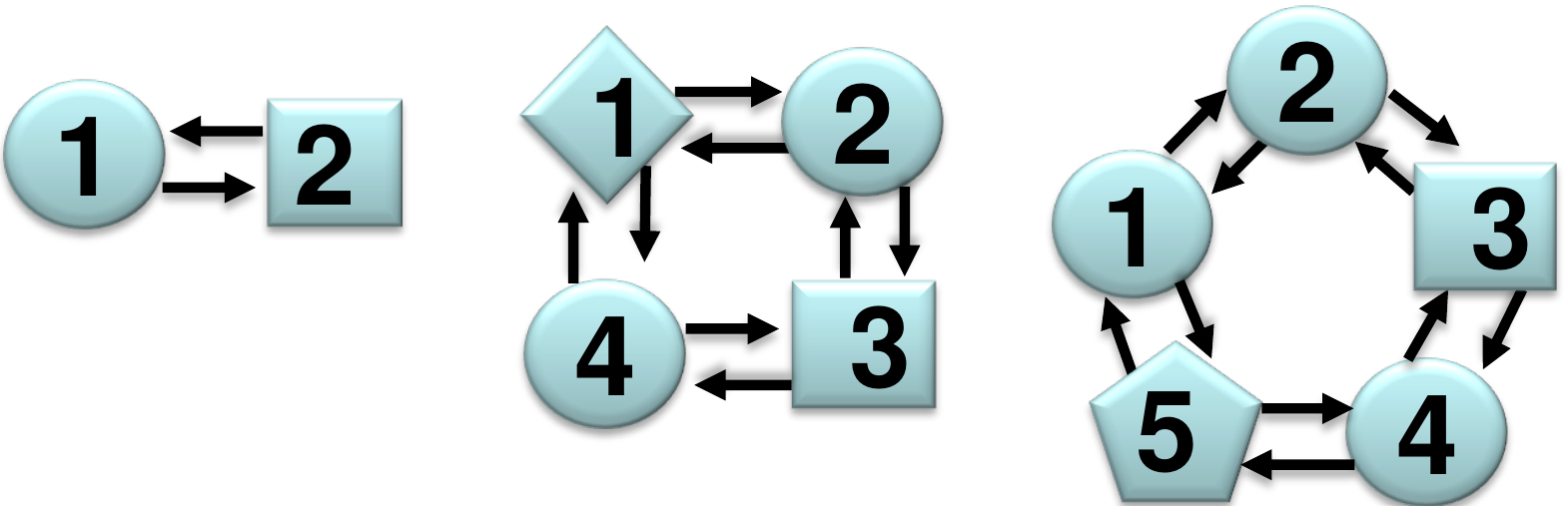}}
 \subfigure[with uni-directional couplings]{\includegraphics{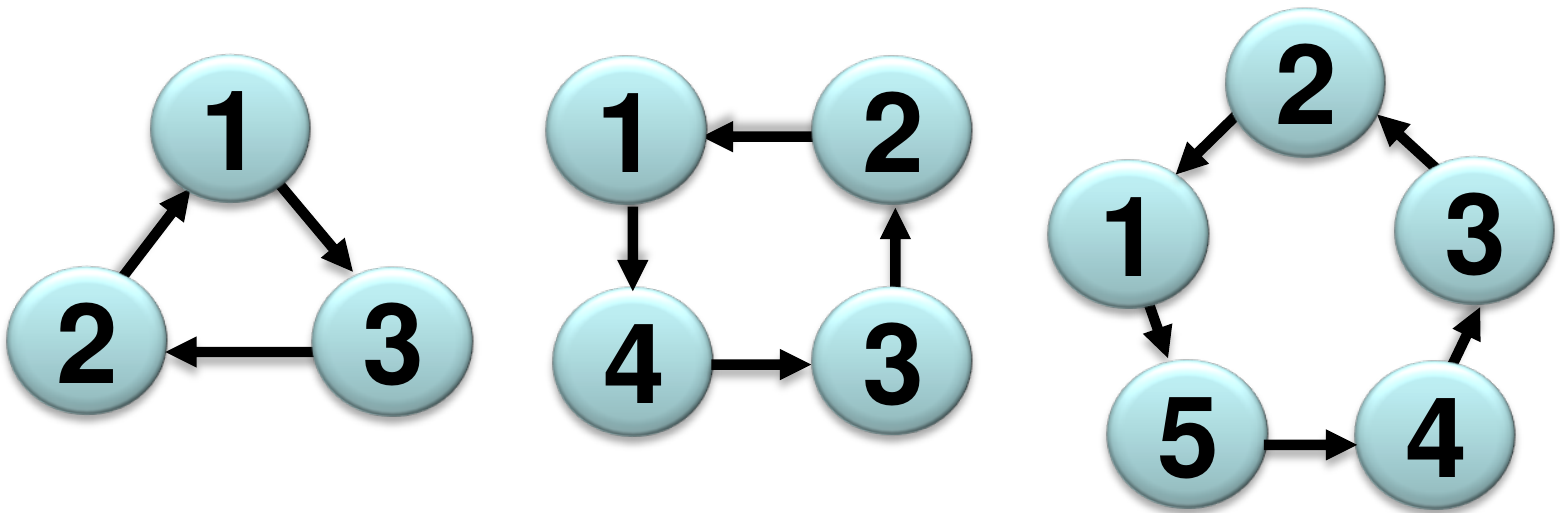}}
 \subfigure[with uni-directional and bi-directional couplings]{\includegraphics{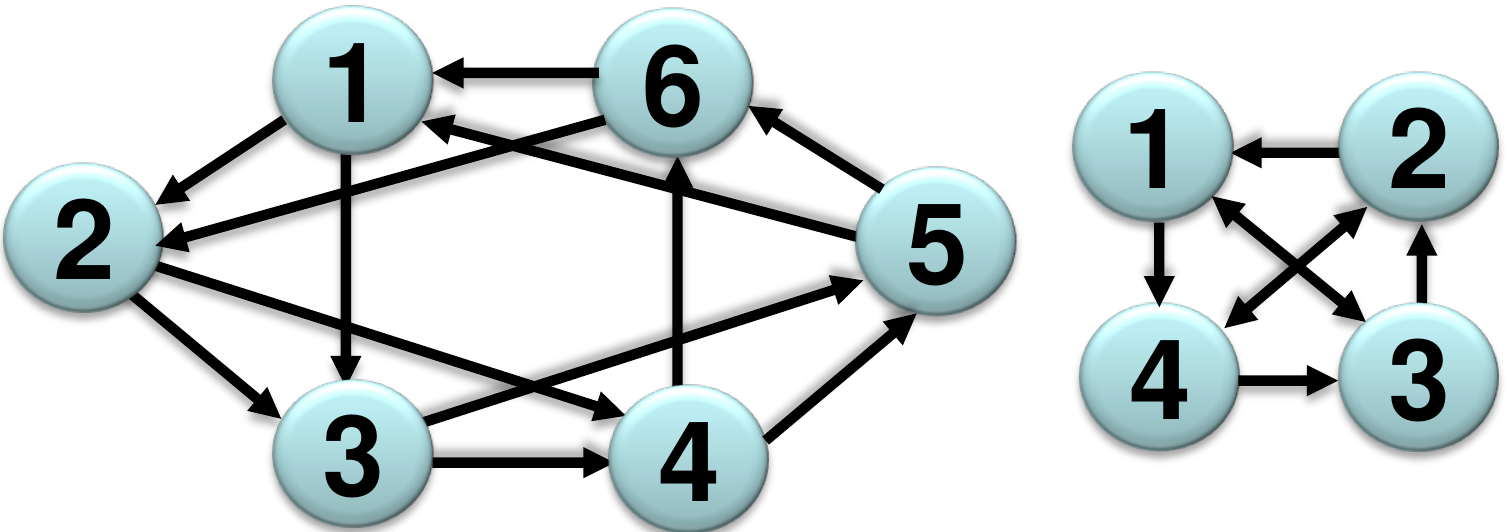}}
 \subfigure[unbalanced graphs with feedback hierarchies]{\includegraphics[width=2.2in]{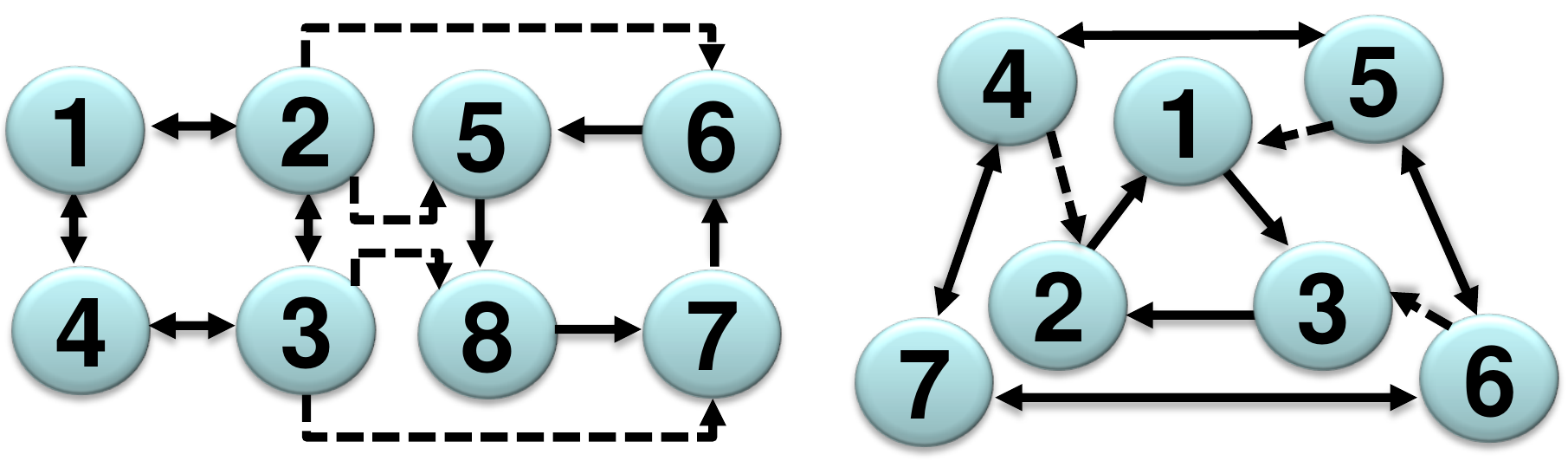}}
  \subfigure[inline configuration]{\includegraphics[width=1.1in]{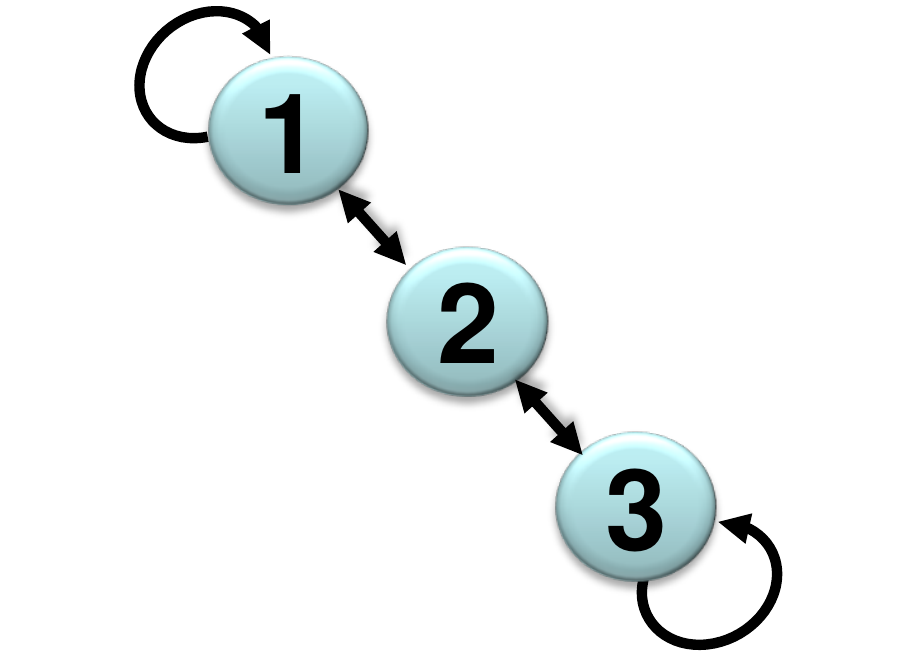}}
 \end{subfigmatrix}
 \caption{Network structures permitted in this paper. Networks in (a-d) are on balanced graphs, and each element has the same number of neighbors (i.e, regular graph). More complex structures, such as the unbalanced graphs shown in (e), can be constructed by concurrent synchronization presented in Section~\ref{sec:concurrent}. The solid lines indicate the local couplings, whereas the dash lines indicate the reference input commands.}\label{robot_diagram}
\end{figure}

Let us first consider the robot networks shown in Fig.~\ref{robot_diagram}(a-d). The following control law, adopted from the single robot control law in \cite{Ref:Slotine}, is proposed for the $i$-th robot in a network comprised of $p$ robots ($p\geq 3$).
\begin{align}\label{tracking_controllerCH5_general}
  \mathbf{\tau}_i=&{\mathbf{M}} (\mathbf{q}_i){\mathbf{\ddot{q}}_{i,r}}+
       \mathbf{C}(\mathbf{q}_i,\mathbf{\dot{q}}_i){\mathbf{\dot{q}}_{i,r}}+\mathbf{g}(\mathbf{q}_i)\\
      &-\mathbf{K}_1(t)(\mathbf{\dot{q}}_i-\mathbf{\dot{q}}_{i,r})+\sum^m_{j\in\mathcal{N}_i(t)}\frac{2}{m}\mathbf{K}_2(t)(\mathbf{\dot{q}}_{j}-\mathbf{\dot{q}}_{j,r})\nonumber
\end{align}
which also permits the uni-directional couplings shown in Figures~\ref{robot_diagram} (c) and (d). Note that $m$, which is the same for each robot, is the number of its neighbors that send the coupling signals to the $i$-th member, and the constant or time-varying set $\mathcal{N}_i(t)$ consists of such neighbors. This generalized form will be discussed again in Section~\ref{sec:concurrent}.

Equation (\ref{tracking_controllerCH5_general}) reduces to the following tracking control law for a
two-way-ring symmetric structure, as shown in Figure~\ref{robot_diagram}(a)
\begin{align}\label{tracking_controllerCH5}
  \mathbf{\tau}_i=&{\mathbf{M}} (\mathbf{q}_i){\mathbf{\ddot{q}}_{i,r}}+
       \mathbf{C}(\mathbf{q}_i,\mathbf{\dot{q}}_i){\mathbf{\dot{q}}_{i,r}}+\mathbf{g}(\mathbf{q}_i)-\mathbf{K}_1(t)(\mathbf{\dot{q}}_i-\mathbf{\dot{q}}_{i,r})   \nonumber \\
       &+\mathbf{K}_2(t)(\mathbf{\dot{q}}_{i-1}-\mathbf{\dot{q}}_{i-1,r})+\mathbf{K}_2(t)(\mathbf{\dot{q}}_{i+1}-\mathbf{\dot{q}}_{i+1,r})
\end{align}
where a uniformly positive definite matrix $\mathbf{K}_1(t)\in
\mathbb{R}^{n\times n}$ is a feedback gain for the $i$-th robot, and
another uniformly positive definite matrix $\mathbf{K}_2(t)\in
\mathbb{R}^{n\times n}$ is a coupling gain with the adjacent members
($i-1$ and $i+1$). While (\ref{tracking_controllerCH5}) constructs a closed graph (e.g., ring), an inline configuration is also permitted since the first and last robot can simply connect the $\mathbf{K}_2(t)$ gain from itself (see Fig.~\ref{robot_diagram}(f) and Section~\ref{sec:inline}). The above control law can also be applied to a
network consisting of $p$ non-identical robots
(Figure~\ref{robot_diagram}(b)), as shall be seen in Section
\ref{sec:nonid_sync}.

While the common desired time-varying trajectory (or the virtual leader dynamics) is denoted by $\mathbf{q}_d(t)$, the reference velocity vector, $\mathbf{\dot{q}}_{i,r}$ is given by
shifting the common desired velocity $\mathbf{\dot{q}}_{d}$ with the
position error:
\begin{equation}\label{ref_vel_vector}
\mathbf{\dot{q}}_{i,r}=\mathbf{\dot{q}}_{d}-\mathbf{\Lambda}\mathbf{\widetilde{q}}_i=\mathbf{\dot{q}}_{d}-\mathbf{\Lambda}(\mathbf{q}_i-\mathbf{q}_{d})
\end{equation}
where $\mathbf{\Lambda}$ is a positive diagonal matrix.

\noindent In contrast with \cite{Ref:robot_sync}, the proposed control law
requires only the local coupling feedback of the most adjacent robots
($i-1$ and $i+1$) for exponential convergence (see Figure
\ref{robot_diagram}). Note that the last ($p$-th) robot is connected
with the first robot to form a ring network as suggested in
\cite{Ref:contraction3}. Moreover, estimates of $\mathbf{\ddot{q}}$
are no longer required.

The closed-loop dynamics using  (\ref{NL_single_compactCH5}) and
(\ref{tracking_controllerCH5}) become
\begin{equation}\label{tracking_closed_simple2}
{\mathbf{M}} (\mathbf{q}_i){\mathbf{\dot{s}}_{i}}+
       \mathbf{C}(\mathbf{q}_i,\mathbf{\dot{q}}_i)\mathbf{s}_{i}+\mathbf{K}_1\mathbf{s}_i-\mathbf{K}_2\mathbf{s}_{i-1}-\mathbf{K}_2\mathbf{s}_{i+1}=\mathbf{0}
\end{equation}
where $\mathbf{s}_i$ denotes the composite variable
$\mathbf{s}_i=\mathbf{\dot q}_i - \mathbf{\dot{q}}_{i,r}$.
\subsection{Modified Laplacian}
Let us define the following $p\times p$ block square matrices:
\begin{equation*}
[\mathbf{L}^{p}_{\mathbf{A},\mathbf{B}}]=\left[\begin{smallmatrix}
\mathbf{A} & \mathbf{B} & \mathbf{0}  & \cdots & \mathbf{B}\\
\mathbf{B} & \mathbf{A} & \mathbf{B} &  \cdots & \mathbf{0}\\
\vdots & \ddots & \ddots &   & \vdots \\
\mathbf{0} & &  \mathbf{B} & \mathbf{A} & \mathbf{B}  \\
\mathbf{B} & \cdots  & \mathbf{0} & \mathbf{B} &
\mathbf{A}
\end{smallmatrix}\right]_{p\times p},  [\mathbf{U}^{p}_{\mathbf{A}}]=\left[\begin{smallmatrix}\mathbf{A}& \mathbf{A} &
\cdots &\mathbf{A}
 \\
 \mathbf{A} & \mathbf{A} & \cdots & \mathbf{A}\\
 \vdots &  \vdots & \ddots & \vdots\\ \mathbf{A}& \mathbf{A} & \cdots & \mathbf{A}
 \end{smallmatrix}\right]_{p\times p}
 \end{equation*}
For a ring structure, defined from the controller in
(\ref{tracking_controllerCH5}),
$[\mathbf{L}^{p}_{\mathbf{A},\mathbf{B}}]$ has only three nonzero
matrix elements in each row (i.e.,
$\mathbf{A},\mathbf{B},\mathbf{B}$). Then, we can write the closed-loop dynamics in
(\ref{tracking_closed_simple2}) in the following block matrix form
\begin{equation}\label{virtual_DN_block}
[\mathbf{M}]\mathbf{\dot x}+
       [\mathbf{C}]\mathbf{x}+\left([\mathbf{L}^{p}_{\mathbf{K}_1,-\mathbf{K}_2}]+[\mathbf{U}^{p}_{\mathbf{K}_2}]\right)\mathbf{x}=[\mathbf{U}^{p}_{\mathbf{K}_2}]\mathbf{x}
\end{equation}
where $[\mathbf{M}]=\text{diag}\left({\mathbf{M}}
(\mathbf{q}_1),\cdots ,  {\mathbf{M}} (\mathbf{q}_p)\right)$,\\
$[\mathbf{C}]=\text{diag}\left(\mathbf{C}(\mathbf{q}_1,\mathbf{\dot{q}}_1),
\cdots ,\mathbf{C}(\mathbf{q}_p,\mathbf{\dot{q}}_p)
 \right)$, $\mathbf{x}=\left(\mathbf{ s}_1^T, \cdots,
\mathbf{s}_p^T\right)^T$.
\begin{definition}$[\mathbf{L}^{p}_{\mathbf{K}_1,-\mathbf{K}_2}]$ can be viewed as the modified Laplacian of the network in the context of graph theory. In
other words, $[\mathbf{L}^{p}_{\mathbf{K}_1,-\mathbf{K}_2}]$
indicates the connectivity with adjacent systems as well as the
strength of the coupling by $\mathbf{K}_2$. Note that $[\mathbf{L}^{p}_{\mathbf{K}_1,-\mathbf{K}_2}]$ can be time-varying due to time-varying $\mathbf{K}_1$ and $\mathbf{K}_2$, or due to the switching topology, $\mathcal{N}_i(t)$ in (\ref{tracking_controllerCH5_general}). \end{definition}
\begin{remark}The network graphs illustrated in Figure
\ref{robot_diagram}(a-d) are \emph{balanced} since the in-degree of each node is equal to the out-degree~\cite{Ref:Murray}. The additional requirement for the stability analysis in this section is that the robots should be on a \emph{regular} graph, where each member has the same number of neighbors. While this paper permits popular local coupling configurations~\cite{Ref:FormationFlying_CHUNG,Jadbabaie} from regular graphs, we will show that the assumption of a regular balanced graph can be relaxed in Section~\ref{sec:concurrent}. In particular, \emph{unbalanced or non-regular graphs due to feedback hierarchies}, as shown in Figure~\ref{robot_diagram}(e), can be employed by concurrent synchronization discussed in Section~\ref{sec:concurrent}.\end{remark}
\begin{remark}It should be noted that the
matrix $[\mathbf{L}^{p}_{\mathbf{K}_1,-\mathbf{K}_2}]$ is different
from the standard Laplacian found in \cite{Ref:Murray}. By
definition, every row sum of the Laplacian matrix on a balanced graph is zero. Hence,
such Laplacian matrix always has a zero eigenvalue corresponding to a
right eigenvector, $\mathbf{1}=(1,1,\cdots,1)^T$.
In contrast, a strictly positive definite
$[\mathbf{L}^{p}_{\mathbf{K}_1,-\mathbf{K}_2}]$ is required for
exponential tracking convergence for the proposed control law in this paper.
In other words, unless otherwise noted,
$[\mathbf{L}^{p}_{\mathbf{K}_1,-\mathbf{K}_2}]$ is assumed to have
no zero eigenvalue. For example, the block matrix for $p=4$ becomes
\begin{equation}[\mathbf{L}^{4}_{\mathbf{K}_1,-\mathbf{K}_2}]=\left[\begin{smallmatrix}+\mathbf{K}_1 & -\mathbf{K}_2 & \mathbf{0} & -\mathbf{K}_2 \\
-\mathbf{K}_2 & +\mathbf{K}_1  & -\mathbf{K}_2  & \mathbf{0} \\
\mathbf{0} & -\mathbf{K}_2 & +\mathbf{K}_1 & -\mathbf{K}_2
\\-\mathbf{K}_2 & \mathbf{0} & -\mathbf{K}_2 & +\mathbf{K}_1
\end{smallmatrix}\right]\end{equation}which is positive definite for
$\mathbf{K}_1-2\mathbf{K}_2>0$. This condition is also true $\forall p,\ p\geq3$.\end{remark}
\subsection{Tracking Stability Analysis}
The following condition should be true for exponential
tracking convergence to the common desired trajectory $\mathbf{q}_d(t)$.
\begin{theorem}\label{Thm:tracking_sync_conv}
\emph{If $[\mathbf{L}^{p}_{\mathbf{K}_1,-\mathbf{K}_2}]$ is uniformly positive
definite:
\begin{equation}
[\mathbf{L}^{p}_{\mathbf{K}_1,-\mathbf{K}_2}]>0, \ \ \forall t
\end{equation}
then every robot follows the desired
trajectory $\mathbf{q}_d(t)$ exponentially fast from any initial
condition.}
%
\end{theorem}
\begin{proof}
We can cancel out the $[\mathbf{U}^{p}_{\mathbf{K}_2}]$ matrix term
in (\ref{virtual_DN_block}) to obtain
\begin{equation}\label{tracking_form}
[\mathbf{M}]\mathbf{\dot x}+
       [\mathbf{C}]\mathbf{x}+[\mathbf{L}^{p}_{\mathbf{K}_1,-\mathbf{K}_2}]\mathbf{x}=\mathbf{0}.
\end{equation}
Equation (\ref{tracking_form}) corresponds to a conventional
tracking problem with a block diagonal matrix of nonlinearly coupled inertia matrices, $[\mathbf{M}]$. We use contraction theory to
prove that $\mathbf{x}$ tends to zero exponentially and globally with
$[\mathbf{L}^{p}_{\mathbf{K}_1,-\mathbf{K}_2}]>0$. Consider the virtual system of $\mathbf{y}$ obtained by replacing
$\mathbf{x}$ with $\mathbf{y}$ in (\ref{tracking_form}).
\begin{equation}
[\mathbf{M}]\mathbf{\dot y}+
       [\mathbf{C}]\mathbf{y}+[\mathbf{L}^{p}_{\mathbf{K}_1,-\mathbf{K}_2}]\mathbf{y}=\mathbf{0}
\end{equation}

This virtual $\mathbf{y}$ system has two particular solutions:
$\mathbf{x}=(\mathbf{ s}_1^T,\cdots,\mathbf{s}_p^T)^T$ and $\mathbf{0}$.
The squared-length analysis with respect to the positive-definite
metric $[\mathbf{M}]$ yields
\begin{align}\label{tracking_proof_squared_length}
&\frac{d}{dt}\left(\delta
\mathbf{y}^T[\mathbf{M}]\delta\mathbf{y}\right)=2\delta
\mathbf{y}^T[\mathbf{M}]\delta\mathbf{\dot y}+\delta
\mathbf{y}^T[\mathbf{\dot M}]\delta\mathbf{y} \\&=-2\delta
\mathbf{y}^T\bigl([\mathbf{C}]\delta\mathbf{y}+[\mathbf{L}^{p}_{\mathbf{K}_1,-\mathbf{K}_2}]\delta\mathbf{y}\bigr)+\delta
\mathbf{y}^T[\mathbf{\dot
M}]\delta\mathbf{y}\nonumber \\&=-2\delta\mathbf{y}^T[\mathbf{L}^{p}_{\mathbf{K}_1,-\mathbf{K}_2}]
\delta\mathbf{y}\nonumber
\end{align}where we used the skew-symmetric property of $[\mathbf{\dot
M}]-2[\mathbf{C}]$.

Accordingly, $[\mathbf{L}^{p}_{\mathbf{K}_1,-\mathbf{K}_2}]>0$ will
make the system contracting ($\delta \mathbf{y}\rightarrow 0$), thus all solutions of $\mathbf{y}$
converge to a single trajectory globally and exponentially fast (Theorems \ref{Thm:contraction} and \ref{Thm:partial}). This in turn
indicates that the composite variable of each robot tends to zero
exponentially ($\mathbf{s}_i\rightarrow \mathbf{0}$). By
definition of
$\mathbf{s}_i=\mathbf{\dot{q}}_i-\mathbf{\dot{q}}_{d}+\mathbf{\Lambda}(\mathbf{q}_i-\mathbf{q}_{d})$,
this show global exponential convergence of $\mathbf{q}_i$ to the common
reference trajectory $\mathbf{q}_{d}(t)$ (see also Theorem \ref{Thm:hierc}).
\end{proof}

Given $\mathbf{K}_1>0,\mathbf{K}_2>0$, it can be shown that a
sufficient condition for the positive-definiteness of
$[\mathbf{L}^{p}_{\mathbf{K}_1,-\mathbf{K}_2}]$ is
$\mathbf{K}_1-\mathbf{K}_2>0$ for $p=2$, and
$\mathbf{K}_1-2\mathbf{K}_2$ for $p \geq 3$.

The next question to be addressed is how to guarantee the
synchronization of the individual dynamics.
\section{Synchronization With/Without Tracking}\label{sec:synsync}
We prove the exponential synchronization
of multiple Lagrangian systems in this section. First, we describe
the difficulties inherent in proving the synchronization of Lagrangian systems in Section \ref{sec:challenges}. We then present
the synchronization proof in Sections \ref{mainproof}. We also show that our method is more general than prior work by reducing our control law to the standard synchronization problem without trajectory tracking in Section \ref{sec:sync_no_tracking}. The adaptive synchronization is presented in Section~\ref{sec:adaptive_sync}.
\subsection{Challenges with Nonlinear Inertia
Matrix}\label{sec:challenges} The difficulties associated with
nonlinear time-varying inertia matrices can be easily demonstrated
with the following two-robot example. The closed-loop dynamics of
two identical robots from (\ref{tracking_closed_simple2}) becomes
\begin{equation}\label{closed_single_tracking}
\begin{split}
  {\mathbf{M}} (\mathbf{q}_1)\mathbf{\dot s}_1+
       \mathbf{C}(\mathbf{q}_1,\mathbf{\dot{q}}_1)\mathbf{s}_1+(\mathbf{K}_1+\mathbf{K}_2)\mathbf{s}_1=\mathbf{u}(t)\\
{\mathbf{M}} (\mathbf{q}_2)\mathbf{\dot s}_2+
       \mathbf{C}(\mathbf{q}_2,\mathbf{\dot{q}}_2)\mathbf{s}_2+(\mathbf{K}_1+\mathbf{K}_2)\mathbf{s}_2=\mathbf{u}(t)\\
\end{split}
\end{equation}
Note that $\mathbf{u}(t)=\mathbf{K}_2(\mathbf{s}_1+{\mathbf{s}_{2}})$ and  that $\mathbf{s}_i$, $i=1,2$ is the composite
variable defined in (\ref{tracking_closed_simple2}).

Direct application of synchronization ({Theorem \ref{Thm:sync}})
appears elusive since we have to prove that (\ref{closed_single_tracking}) are
contracting in the same metric while preserving the input
symmetry~\cite{Ref:contraction_sync}. Hence, this can be viewed as a higher order contraction problem~\cite{Ref:Slotine_highorder}. For example, multiplying
(\ref{closed_single_tracking}) by $\mathbf{M}^{-1}$ breaks the input
symmetry: i.e., ${\mathbf{M}}^{-1}(\mathbf{q}_1)\mathbf{u}(t)\neq
{\mathbf{M}}^{-1}(\mathbf{q}_2)\mathbf{u}(t)$. In essence,
${\mathbf{M}}(\mathbf{q}_1)\neq {\mathbf{M}}(\mathbf{q}_2)$ makes
this problem intractable in general.

Instead, assume that ${\mathbf{M}}(\mathbf{q})$ becomes a constant matrix, thereby making $\mathbf{C}(\mathbf{q},\mathbf{\dot{q}})$
zero. Then, we can easily prove $\mathbf{s}_1$ and $\mathbf{s}_2$
tend to each other from
\begin{equation}\label{prior_work_examp}
  {\mathbf{M}} \mathbf{\dot s}_1+
       (\mathbf{K}_1+\mathbf{K}_2)\mathbf{s}_1=\mathbf{u}(t),\ \
{\mathbf{M}} \mathbf{\dot s}_2+
       (\mathbf{K}_1+\mathbf{K}_2)\mathbf{s}_2=\mathbf{u}(t)
\end{equation}

\noindent The following virtual $\mathbf{y}$-system with the common input
$\mathbf{u}(t)$
\begin{equation}
{\mathbf{M}} \mathbf{\dot y}+
       (\mathbf{K}_1+\mathbf{K}_2)\mathbf{y}=\mathbf{u}(t)
\end{equation}
is partially contracting with $\mathbf{K}_1+\mathbf{K}_2>0$ (see Theorem~\ref{Thm:partial}). Hence, its
particular solutions $\mathbf{s}_1$ and $\mathbf{s}_2$ tend to each
other exponentially fast according to the synchronization theorem
(Theorem \ref{Thm:sync}). Without loss of generality, this result
can easily be extended to arbitrarily large networks. The
synchronization of a large network with a constant metric, as seen in (\ref{prior_work_examp}), is already
discussed in \cite{Ref:contraction3} using contraction analysis.

We now turn to a more difficult problem focused on the
synchronization of two robots with non-constant nonlinear metrics
(${\mathbf{M}}(\mathbf{q}_1)\neq {\mathbf{M}}(\mathbf{q}_2)$).
\subsection{Contraction with Multiple Time Scales}
In this section, we show that we can render the system synchronized first, then follow
the common trajectory by tuning the gains properly. This indicates that there exist two different time scales in the
closed-loop systems constructed with the proposed controllers. This multi-time-scale behavior will be
exploited in the subsequent sections.

Recall the closed-loop dynamics given in (\ref{tracking_form}). Since $[\mathbf{L}^{p}_{\mathbf{K}_1,-\mathbf{K}_2}]$ is a real
symmetric matrix, we can perform the spectral decomposition~\cite{Ref:SpectralGraph}. This is a special case of the concurrent synchronization~\cite{Ref:contraction_sync} that corresponds to convergence to a flow invariant subspace (the eigenspace).
\begin{equation}\label{L_diagonalization}
[\mathbf{L}^{p}_{\mathbf{K}_1,-\mathbf{K}_2}]=\mathbf{V}[\mathbf{D}]\mathbf{V}^T
\end{equation}
where $[\mathbf{D}]$ is a block diagonal matrix, and the square matrix $\mathbf{V}$ is composed of the orthonormal eigenvectors such that $\mathbf{V}^T\mathbf{V}=\mathbf{V}\mathbf{V}^T=\mathbf{I}_{pn}$, since the symmetry
of $[\mathbf{L}^{p}_{\mathbf{K}_1,-\mathbf{K}_2}]$ gives rise to real eigenvalues and orthonormal
eigenvectors~\cite{Ref:Strang}.

Pre-multiplying (\ref{tracking_form}) by $\mathbf{V}^T$ and
setting $\mathbf{V}^T\mathbf{x}=\mathbf{z}$ result in
\begin{equation}\label{premult2}
\left(\mathbf{V}^T[\mathbf{M}]\mathbf{V}\right)\mathbf{\dot z}+
       \left(\mathbf{V}^T[\mathbf{C}]\mathbf{V}\right)\mathbf{
z}+[\mathbf{D}]\mathbf{z}=\mathbf{0}
\end{equation}

Then, we can develop the squared-length analysis, as in (\ref{tracking_proof_squared_length}). This follows from the fact that
$\left(\mathbf{V}^T[\mathbf{M}]\mathbf{V}\right)$ is always
symmetric positive definite due to a symmetric positive definite $[\mathbf{M}]$.

Since the modified Laplacian $[\mathbf{L}^{p}_{\mathbf{K}_1,-\mathbf{K}_2}]$ represents a \emph{regular} graph, where each member has the same number of neighbors ($=2$ for $p\geq3$),
\begin{equation}
[\mathbf{1}]=\frac{1}{\sqrt{p}}[\mathbf{I}_n,\mathbf{I}_n,\cdots,\mathbf{I}_n]^T
\end{equation}
is the $pn\times n$ block column matrix of eigenvectors associated with the eigenvalues $\lambda(\mathbf{K}_1-2\mathbf{K}_2)$ for $p\geq3$. Note that $\mathbf{I}_n$ denotes the $n \times n$ identity matrix, and the $[\mathbf{1}]$ matrix consists of $p$ matrices of $\mathbf{I}_n$. The eigenvector matrix $[\mathbf{1}]$ represents the common reference trajectory tracking state.

We can define a $pn \times (p-1)n$ matrix $\mathbf{V}_{sync}$ that consists of the orthonormal eigenvectors other than $[\mathbf{1}]$ such that
\begin{align}\label{def_Vsync}
&\mathbf{V}^T\mathbf{V}=\begin{pmatrix}[\mathbf{1}]^T\\\mathbf{V}_{sync}^T\end{pmatrix}\begin{pmatrix}[\mathbf{1}] & \mathbf{V}_{sync}\end{pmatrix}\\&=\begin{bmatrix}[\mathbf{1}]^T[\mathbf{1}]&[\mathbf{1}]^T\mathbf{V}_{sync}\\\mathbf{V}_{sync}^T[\mathbf{1}]&\mathbf{V}_{sync}^T\mathbf{V}_{sync}\end{bmatrix}=
\begin{bmatrix}\mathbf{I}_n& \mathbf{0}_{n \times (p-1)n}  \\\mathbf{0}_{(p-1)n \times n} & \mathbf{I}_{(p-1)n}\end{bmatrix}\nonumber
\end{align}
where we used the orthogonality between $[\mathbf{1}]$ and $\mathbf{V}_{sync}$.

Hence, the block diagonal matrix $[\mathbf{D}]$, which represents the eigenvalues of $[\mathbf{L}^{p}_{\mathbf{K}_1,-\mathbf{K}_2}]$, can be partitioned from (\ref{L_diagonalization})
\begin{align}
[\mathbf{D}]&=\mathbf{V}^T[\mathbf{L}^{p}_{\mathbf{K}_1,-\mathbf{K}_2}]\mathbf{V}\nonumber \\
&=\begin{bmatrix}[\mathbf{1}]^T[\mathbf{L}^{p}_{\mathbf{K}_1,-\mathbf{K}_2}][\mathbf{1}]&[\mathbf{1}]^T[\mathbf{L}^{p}_{\mathbf{K}_1,-\mathbf{K}_2}]\mathbf{V}_{sync}\\ \mathbf{V}_{sync}^T[\mathbf{L}^{p}_{\mathbf{K}_1,-\mathbf{K}_2}][\mathbf{1}]&\mathbf{V}_{sync}^T[\mathbf{L}^{p}_{\mathbf{K}_1,-\mathbf{K}_2}]\mathbf{V}_{sync}\end{bmatrix}\nonumber\\
&=\begin{bmatrix}\mathbf{D}_1& \mathbf{0}_{n \times (p-1)n}\\ \mathbf{0}_{(p-1)n \times n} & \mathbf{D}_2\end{bmatrix}
\end{align}

It should be emphasized that $\mathbf{D}_1$, which equals $\mathbf{K}_1-2\mathbf{K}_2$ for $p\geq3$, represents the tracking gain, while $\mathbf{D}_2$ corresponds to the synchronization gain. We can choose the diagonal control gain matrices $\mathbf{K}_1$ and $\mathbf{K}_2$ such that
\begin{equation}
\mathbf{D}_2=\mathbf{V}_{sync}^T[\mathbf{L}^{p}_{\mathbf{K}_1,-\mathbf{K}_2}]\mathbf{V}_{sync}>\mathbf{D}_1=[\mathbf{1}]^T[\mathbf{L}^{p}_{\mathbf{K}_1,-\mathbf{K}_2}][\mathbf{1}],\nonumber
\end{equation}
thereby ensuring that the robots synchronize faster than they follow the common desired trajectory.

\begin{figure}
  \includegraphics[width=3.2in]{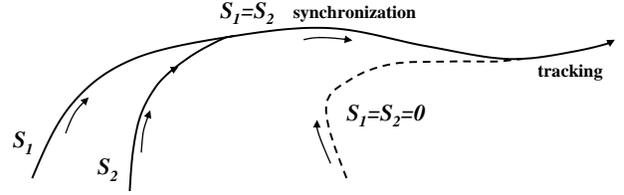}\\
  \caption{Multiple timescales of synchronization (faster) and tracking (slower). The dashed line indicates the desired trajectory, and arrows indicate increasing time. The drawing is conceptual, since strictly speaking $\mathbf{s}_1$ and $\mathbf{s}_2$ synchronize exponentially.}\label{sync_track_graph}
\end{figure}
\noindent This multi timescale behavior is graphically illustrated in
Figure~\ref{sync_track_graph}. The figure depicts that
$\mathbf{s}_1$ and $\mathbf{s}_2$ synchronize first, then they
converge to the desired trajectory while staying together. This
observation motivates separation of the two different time scales,
namely $\mathbf{D}_1$ for tracking and $\mathbf{D}_2$ for synchronization.
\subsection{Stability Analysis of Exponential Synchronization}\label{mainproof}
Using the results from the previous sections, we present the main theorem on synchronization.
\begin{theorem}\label{Thm:tracking_sync}
\emph{Assume that the conditions in Theorem~\ref{Thm:tracking_sync_conv} are
true, so that the individual dynamics are exponentially tracking the
common desired trajectory. A swarm of $p$ robots synchronizes
exponentially from any initial conditions if $\exists $ diagonal
matrices $\mathbf{K}_1>0$, $\mathbf{K}_2>0$, $\mathbf{\Lambda}>0$ such that}
\begin{equation}\label{sync_condition}
\mathbf{D}_2=\mathbf{V}_{sync}^T[\mathbf{L}^{p}_{\mathbf{K}_1,-\mathbf{K}_2}]\mathbf{V}_{sync}>0, \ \ \forall t
\end{equation}
\end{theorem}
\begin{proof}
Consider the virtual system of $\mathbf{y}$ from (\ref{premult2}):
\begin{align}\label{virtual_chap5}
\left(\mathbf{V}^T[\mathbf{M}]\mathbf{V}\right)\mathbf{\dot y}+
       \left(\mathbf{V}^T[\mathbf{C}]\mathbf{V}\right)\mathbf{
y}+[\mathbf{D}]\mathbf{y}=\mathbf{0},
\end{align}
This has $\mathbf{y}=\mathbf{V}^T\mathbf{x}$ and $\mathbf{y}=\mathbf{0}$ as particular solutions, which can be written in terms of $\mathbf{y}=\left(\mathbf{y}_t^T,\mathbf{y}_s^T\right)^T$:
\begin{equation}\label{proof_particular_solution}
\begin{pmatrix}\mathbf{y}_t=[\mathbf{1}]^T\mathbf{x}\\ \mathbf{y}_s=\mathbf{V}_{sync}^T\mathbf{x}\end{pmatrix}\ \ \ \text{and} \ \ \ \begin{pmatrix}\mathbf{y}_t=\mathbf{0}\\ \mathbf{y}_s=\mathbf{0}\end{pmatrix}\end{equation}

For $[\mathbf{D}]>0$, we
can show that the above virtual system is contracting (i.e., $\delta \mathbf{y}\rightarrow \mathbf{0}$ globally and exponentially). We take the
symmetric positive definite block matrix
$\mathbf{V}^T[\mathbf{M}]\mathbf{V}$
as our contraction metric.

Performing the squared-length analysis with respect to this metric
yields
\begin{align}\label{virtual_length_good}
\frac{d}{dt}\delta
\mathbf{y}^T&\left(\mathbf{V}^T[\mathbf{M}]\mathbf{V}\right)\delta\mathbf{y} =-2\delta
\mathbf{y}^T\Bigl( \left(\mathbf{V}^T[\mathbf{C}]\mathbf{V}\right)\delta\mathbf{y}+[\mathbf{D}]\delta\mathbf{y}\Bigr)\nonumber \\&+\delta
\mathbf{y}^T\left(\mathbf{V}^T[\mathbf{\dot M}]\mathbf{V}\right)\delta\mathbf{y}=
-2\delta \mathbf{y}^T[\mathbf{D}]\delta\mathbf{y}
\end{align}
where we used the skew-symmetric property of $\left(\mathbf{V}^T[\mathbf{\dot M}]\mathbf{V}\right)-2 \left(\mathbf{V}^T[\mathbf{C}]\mathbf{V}\right)$.

The above equation can be rewritten in terms of two different time scales
\begin{align}\label{two_different_scales}
&\frac{d}{dt}\begin{pmatrix}\delta\mathbf{y}_t \\ \delta\mathbf{y}_s\end{pmatrix}^T\begin{bmatrix}[\mathbf{1}]^T[\mathbf{M}][\mathbf{1}]&[\mathbf{1}]^T[\mathbf{M}]\mathbf{V}_{sync}\\ \mathbf{V}_{sync}^T[\mathbf{M}][\mathbf{1}]&\mathbf{V}_{sync}^T[\mathbf{M}]\mathbf{V}_{sync}\end{bmatrix}\begin{pmatrix}\delta\mathbf{y}_t \\ \delta\mathbf{y}_s\end{pmatrix}\\&=-2\begin{pmatrix}\delta\mathbf{y}_t \\ \delta\mathbf{y}_s\end{pmatrix}^T\begin{bmatrix}\mathbf{D}_1& \mathbf{0}\\ \mathbf{0} & \mathbf{D}_2\end{bmatrix}\begin{pmatrix}\delta\mathbf{y}_t \\ \delta\mathbf{y}_s\end{pmatrix}\nonumber
\end{align}

If $\mathbf{D}_1>0$ and $\mathbf{D}_2>0$, the combined virtual system in (\ref{virtual_chap5})
is contracting. In other words, $\delta \mathbf{y} \rightarrow \mathbf{0}$ exponentially fast. This in
turn implies that all solutions of $\mathbf{y}$ tend
to the single trajectory. In particular, the tracking ($\delta \mathbf{y}_t \rightarrow \mathbf{0}$) is associated with $\mathbf{D}_1$, and synchronization ($\delta \mathbf{y}_s \rightarrow \mathbf{0}$) is associated with $\mathbf{D}_2$.

As a result,
$[\mathbf{1}]^T\mathbf{x}=1/\sqrt{p}(\mathbf{s}_1+\cdots+\mathbf{s}_p)$ and
$\mathbf{V}_{sync}^T\mathbf{x}$ from (\ref{proof_particular_solution}) tend to zero exponentially. Note that $\mathbf{s}_1,\cdots,\mathbf{s}_p \rightarrow \mathbf{0}$ has already been proven with $\mathbf{D}_1>0$ for Theorem \ref{Thm:tracking_sync_conv}, which is a sufficient condition to make the sum of the composite variables also tend to zero (i.e., $[\mathbf{1}]^T\mathbf{x}\rightarrow\mathbf{0}$). What is new in this section is the synchronization $\mathbf{V}_{sync}^T\mathbf{x}\rightarrow \mathbf{0}$ and its convergence rate with $\mathbf{D}_2>0$.

It is straightforward to show that $\mathbf{V}_{sync}^T\mathbf{x} \rightarrow \mathbf{0}$ and $\mathbf{\Lambda}>0$ also
hierarchically make $\mathbf{q}_1,\cdots,\mathbf{q}_p$ synchronize
globally exponentially fast (see Theorem \ref{Thm:hierc}). This can be verified by the following contracting dynamics constructed from (\ref{tracking_closed_simple2})
\begin{equation}\label{q_contraction}
\mathbf{V}_{sync}^T\{\mathbf{\dot q}\}+\left(\mathbf{V}_{sync}^T[\mathbf{\Lambda}]\mathbf{V}_{sync}\right)\mathbf{V}_{sync}^T\{\mathbf{q}\}=\mathbf{V}_{sync}^T\mathbf{x}\rightarrow \mathbf{0}
\end{equation}
where $\{\mathbf{q}\}=(\mathbf{q}_1^T,\cdots,\mathbf{q}_p^T)^T$. Note that the orthonormal vectors $\mathbf{V}_{sync}$ canceled the common input term $\mathbf{\dot q}_d+\mathbf{\Lambda}\mathbf{q}_d$. Also, $[\mathbf{\Lambda}]$ is a block diagonal matrix of $\mathbf{\Lambda}>0$, thereby yielding $\mathbf{V}_{sync}^T[\mathbf{\Lambda}][\mathbf{1}]=\mathbf{0}$ from
\begin{equation}
\mathbf{V}_{sync}\mathbf{V}_{sync}^T\{\mathbf{q}\}+[\mathbf{1}][\mathbf{1}]^T\{\mathbf{q}\}=\{\mathbf{q}\}
\end{equation}

Consequently, $\mathbf{V}_{sync}^T\{\mathbf{q}\}\rightarrow\mathbf{0}$ implies the synchronization of the original state variable $\mathbf{q}_i$,  $i=1,\cdots,p$.
This also implies that the diagonal terms of the metric,
$\mathbf{V}_{sync}^T[\mathbf{M}][\mathbf{1}]$, tend
to zero exponentially, thereby eliminating the coupling of the
inertia term $\mathbf{V}^T[\mathbf{M}]\mathbf{V}$ in (\ref{virtual_length_good}--\ref{two_different_scales}). This completes the proof of Theorem~\ref{Thm:tracking_sync}.\end{proof}

\begin{remark}We assume here that $\mathbf{q}_d(t)$ is identical for each member. If $\mathbf{q}_d(t)$ were different for each dynamics, $\mathbf{s}_i \rightarrow \mathbf{s}_j$ would imply the synchronization of $\mathbf{q}_i-\mathbf{q}_j$ to the difference of the desired trajectories, which would be useful to construct phase synchronization of oscillatory trajectories (e.g., see~\cite{Ref:FormationFlying_CHUNG,Ref:contraction_sync}).\end{remark}

We can use Theorem \ref{Thm:tracking_sync} without finding $\mathbf{V}_{sync}$ as follows.
\begin{corollary}\label{Corollary1}
\emph{The following condition, in lieu of (\ref{sync_condition}), verifies Theorem~\ref{Thm:tracking_sync}.
\begin{equation}
[\mathbf{L}^{p}_{\mathbf{K}_1,-\mathbf{K}_2}]+[\mathbf{U}^{p}_{\mathbf{K}_2}]>0, \ \ \forall t
\end{equation}}
\end{corollary}
\begin{proof}
The block matrix $[\mathbf{U}^{p}_{\mathbf{K}_2}]$ also has $[\mathbf{1}]$ as its eigenvector. We multiply $[\mathbf{K}_{sync}]=[\mathbf{L}^{p}_{\mathbf{K}_1,-\mathbf{K}_2}]+[\mathbf{U}^{p}_{\mathbf{K}_2}]$ by its orthonormal eigenvectors other than $[\mathbf{1}]$:
\begin{equation}
[\mathbf{K}_{sync}]\mathbf{V}_{sync}=[\mathbf{L}^{p}_{\mathbf{K}_1,-\mathbf{K}_2}]\mathbf{V}_{sync}=\mathbf{V}_{sync}\mathbf{D}_2
\end{equation}
which shows that $\mathbf{V}_{sync}$ also represents the orthonormal eigenvectors of $[\mathbf{K}_{sync}]=[\mathbf{L}^{p}_{\mathbf{K}_1,-\mathbf{K}_2}]+[\mathbf{U}^{p}_{\mathbf{K}_2}]$. In other words, $\mathbf{D}_2$ corresponds to the eigenvalues of $[\mathbf{K}_{sync}]$. The remaining eigenvalue of $[\mathbf{K}_{sync}]$ that is associated with $[\mathbf{1}]$ is $\mathbf{K}_1+(p-2)\mathbf{K}_2$, which is greater than the tracking eigenvalue $\mathbf{D}_1=\mathbf{K}_1-2\mathbf{K}_2$ for $p\geq3$. Hence, the synchronization occurs with $[\mathbf{L}^{p}_{\mathbf{K}_1,-\mathbf{K}_2}]+[\mathbf{U}^{2}_{\mathbf{K}_2}]>0$.\end{proof}

Theorem~\ref{Thm:tracking_sync} and Corollary~\ref{Corollary1} correspond to synchronization with stable
tracking. As shall be seen in Section \ref{sec:nonid_sync}, multiple
dynamics need not be identical to achieve stable synchronization. It is useful to note that the above condition corresponds to
$\mathbf{K}_1+\mathbf{K}_2>0$ for two-robot and three-robot networks
($p=2,3$).
\begin{remark}\emph{Robust Synchronization.} By extending \cite{Ref:contraction1}, we can show that coupled contracting dynamics have the property of robustness to bounded disturbances~\cite{Ref:sjchung_PhD}. It should be emphasized that exponential stability of contraction analysis facilitates such a perturbation analysis. In general, the proof of robustness with asymptotic convergence is more involved, or even leads to instability with respect to a certain class of perturbations (see p. 350 in~\cite{Khalil:2002}).\end{remark}

In the next section, we show that a network of
multiple robots can synchronize even without stable tracking.
\subsection{Synchronization Without a Common Reference Trajectory}\label{sec:sync_no_tracking}
We now turn into the more standard synchronization problem where the tracking gain is zero ($\mathbf{D}_1=\mathbf{K}_1-2\mathbf{K}_2=0$ for $p\geq3$), which fails the exponential
tracking stability condition in Theorem \ref{Thm:tracking_sync_conv}. In this case, the modified Laplacian $[\mathbf{L}^{p}_{\mathbf{K}_1,-\mathbf{K}_2}]$ reduces to the standard weighted Laplacian whose row sum is zero. For synchronization to the weighted average of the initial conditions, we do not require the common desired trajectory $\mathbf{q}_d$, and $\mathbf{q}_d$ can simply be set to zero as follows:
\begin{equation}\label{defsemicontractings}
\mathbf{\dot q}_{i,r}=-\mathbf{\mathbf{\Lambda}} \mathbf{q}_i,\ \ \ \ \ \ \mathbf{s}_i=\mathbf{\dot q}_i+\mathbf{\Lambda}\mathbf{q}_i
\end{equation}

In other words, our control strategy represents a more generalized framework for the synchronization of multiple Lagrangian systems.
\begin{theorem}\label{Thm:tracking_sync_unstable}
\emph{Suppose that the individual tracking dynamics are \emph{indifferent}. Hence, the conditions in Theorem~\ref{Thm:tracking_sync_conv} are \emph{not} true. Nevertheless, a swarm of $p$
\emph{identical} robots asymptotically synchronize from any
initial conditions if $\exists$ diagonal matrices $\mathbf{K}_1>0$,
$\mathbf{K}_2$ such that}
\begin{equation}\label{semi_contraction_condition}
[\mathbf{L}^{p}_{\mathbf{K}_1,-\mathbf{K}_2}]+[\mathbf{U}^{p}_{\mathbf{K}_2}]>0, \ \ \forall t
\end{equation}
\emph{For indifferent tracking, a common desired trajectory $\mathbf{q}_d(t)$ is no longer required.}
\end{theorem}
\begin{proof}The
combined virtual system per se is then semi-contracting
(see \cite{Ref:contraction1}) since the squared-length
analysis in (\ref{virtual_length_good}-\ref{two_different_scales}) yields the negative
semi-definite matrix:
\begin{align}\label{semicontractingvirtual}
\dot V&=-2\begin{pmatrix}\delta\mathbf{y}_t \\ \delta\mathbf{y}_s\end{pmatrix}^T\begin{bmatrix}\mathbf{0}& \mathbf{0}\\ \mathbf{0} & \mathbf{D}_2\end{bmatrix}\begin{pmatrix}\delta\mathbf{y}_t \\ \delta\mathbf{y}_s\end{pmatrix} \leq0
\end{align}

\noindent While $\delta \mathbf{y}_t$, representing the tracking dynamics,
remains in a finite ball due to $\mathbf{D}_1=\mathbf{0}$,
$\delta \mathbf{y}_s$ tends to zero asymptotically due to
$\mathbf{D}_2>0$. This result can be proven as
follows. Note that $\ddot V=
-4\delta\mathbf{y}_s^T\mathbf{D}_2\delta\mathbf{\dot
y}_s$. $\dot V$ is uniformly continuous since a bounded
$\delta\mathbf{\dot y}_s$ from (\ref{virtual_chap5}) leads to
a bounded $\ddot V$. Due to $\dot V \leq 0$, the use of Barbalat's
lemma~\cite{Ref:Slotine} verifies that $\dot V \rightarrow 0$ as
$t\rightarrow \infty$. This implies that $\delta\mathbf{y}_s$ tends
to zero asymptotically fast. As a result, $\mathbf{V}_{sync}^T\mathbf{x}$ tends to zero asymptotically. From the hierarchical combination discussed in (\ref{q_contraction}), this also implies $\mathbf{V}_{sync}^T\{\mathbf{q}\}\rightarrow 0$. This completes the proof of Theorem~\ref{Thm:tracking_sync_unstable} with indifferent tracking.\end{proof}

This will eventually decouple
the metric matrix with $\mathbf{\Lambda}>0$, as seen in
(\ref{two_different_scales}), since $\mathbf{V}_{sync}^T[\mathbf{M}][\mathbf{1}]$ tends to zero
simultaneously as $\mathbf{V}_{sync}^T\{\mathbf{q}\}\rightarrow \mathbf{0}$. As a
result, when ${\mathbf{M}} (\mathbf{q}_i)-{\mathbf{M}}
(\mathbf{q}_j)$ is sufficiently close to zero, the convergence of
$\delta\mathbf{y}_s \rightarrow 0$ turns exponential.
%
%
\begin{remark}\emph{Fast Inhibition.} The dynamics of a large network with semi-contracting stability as in (\ref{semicontractingvirtual}) can be instantaneously transformed to contracting dynamics by the addition of a single inhibitory coupling link. In other words, a single inhibitory link will also make $\delta \mathbf{y}_t \rightarrow 0$. For instance, we can add a single inhibitory link between two arbitrary elements $a$ and $b$ while we keep the rest of the elements the same:~\cite{Ref:contraction3}
\begin{align}\label{tracking_controller_inhibition}
  \mathbf{\tau}_a=&{\mathbf{M}} (\mathbf{q}_a){\mathbf{\ddot{q}}_{ar}}+
       \mathbf{C}(\mathbf{q}_a,\mathbf{\dot{q}}_a){\mathbf{\dot{q}}_{ar}}+\mathbf{g}(\mathbf{q}_a)\nonumber \\
       &-2\mathbf{K}_2\mathbf{s}_a+\mathbf{K}_2\mathbf{s}_{a-1}+\mathbf{K}_2\mathbf{s}_{a+1}-\mathbf{K}(\mathbf{s}_a+\mathbf{s}_b)\\
\mathbf{\tau}_b=&{\mathbf{M}} (\mathbf{q}_b){\mathbf{\ddot{q}}_{br}}+
       \mathbf{C}(\mathbf{q}_b,\mathbf{\dot{q}}_b){\mathbf{\dot{q}}_{br}}+\mathbf{g}(\mathbf{q}_b)\nonumber \\
       &-2\mathbf{K}_2\mathbf{s}_b+\mathbf{K}_2\mathbf{s}_{b-1}+\mathbf{K}_2\mathbf{s}_{b+1}-\mathbf{K}(\mathbf{s}_a+\mathbf{s}_b)\nonumber
\end{align}
where $2\mathbf{K}_2$ is substituted for $\mathbf{K}_1$, and $\mathbf{K}>0$.

Hence, we can straightforwardly show that $[\mathbf{L}^{p}_{\mathbf{K}_1,-\mathbf{K}_2}]$ is now strictly positive definite, in contrast with the original semi-contracting system. As a result, the closed-loop system is contracting, resulting in $\mathbf{s}_i \rightarrow \mathbf{0}$ and $\mathbf{q}_i \rightarrow \mathbf{0}$ from (\ref{defsemicontractings}). This fast inhibition is useful to rapidly destroy unwanted synchronized oscillation of robots.\end{remark}
\subsection{Adaptive Synchronization}\label{sec:adaptive_sync}
We present the adaptive version of the proposed control law that adapts to the unknown parametric uncertainties of the robot dynamic models. Consider the following adaptive control law~\cite{Jouffroy,Ref:Slotine,Ref:contraction2} that has the same
local coupling structure as the proposed control law in
(\ref{tracking_controllerCH5}):
\begin{align}\label{adaptive_sync_law}
\mathbf{\tau}_i&=\mathbf{\hat{M}}_i(\mathbf{q}_i)\mathbf{\ddot{q}}_{i,r}+
       \mathbf{\hat{C}}_i(\mathbf{q}_i,\mathbf{\dot q}_i)\mathbf{\dot{q}}_{i,r}+\mathbf{\hat{g}}_i(\mathbf{q}_i)
      -\mathbf{K}_1\mathbf{s}_i\\&+\mathbf{K}_2\mathbf{s}_{i-1}+\mathbf{K}_2\mathbf{s}_{i+1}
      =\mathbf{Y}_i\hat{\mathbf{a}}_i
-\mathbf{K}_1\mathbf{s}_i+\mathbf{K}_2\mathbf{s}_{i-1}+\mathbf{K}_2\mathbf{s}_{i+1}\nonumber
\end{align}
where $\mathbf{s}_i=\mathbf{\dot q}_i -
\mathbf{\dot{q}}_{i,r}$. Also, the parameter estimate $\mathbf{\hat a}_i$ for the $i$-th member is
updated by the correlation integral
\begin{equation}\label{adaptive_defnition}
\mathbf{\dot{\hat a}}_i=-\mathbf{\Gamma}\mathbf{Y}_i^T\mathbf{s}_i
\end{equation}
where $\mathbf{\Gamma}$ is a symmetric positive definite matrix.
Hence, the closed-loop system for a network comprised of $p$
non-identical robots can be written as
\begin{align}\label{adaptive_block}
\begin{bmatrix}[\mathbf{M}] & \mathbf{0}\\ \mathbf{0} & [\mathbf{\Gamma^{-1}}] \end{bmatrix}\begin{pmatrix}\mathbf{\dot x}\\\{ \mathbf{\dot{\tilde a}}\}\end{pmatrix}
      &+\begin{bmatrix} [\mathbf{C}]&\mathbf{0} \\ \mathbf{0}&\mathbf{0} \end{bmatrix}\begin{pmatrix}\mathbf{ x}\\ \{\mathbf{{\tilde a}}\}\end{pmatrix}\\
      &+\begin{bmatrix}[\mathbf{L}^{p}_{\mathbf{K}_1,-\mathbf{K}_2}] & -[\mathbf{Y}] \\ [\mathbf{Y}]^T& \mathbf{0}\end{bmatrix}\begin{pmatrix}\mathbf{ x}\\ \{\mathbf{{\tilde a}}\}\end{pmatrix}=\mathbf{0}\nonumber
\end{align}
where $[\mathbf{M}]$ and $[\mathbf{C}]$ are the block diagonal matrices of $\mathbf{M}_i(\mathbf{q}_i)$ and $\mathbf{C}_i(\mathbf{q}_i,\mathbf{\dot q}_i)$, $i=1,\cdots,p$, as defined in (\ref{virtual_DN_block}). The additional block diagonal matrices are defined from (\ref{adaptive_defnition}) such that
$[\mathbf{\Gamma^{-1}}]=\text{diag}(\mathbf{\Gamma}^{-1},\mathbf{\Gamma}^{-1}\cdots,\mathbf{\Gamma}^{-1})_p$, $[\mathbf{Y}]=\text{diag}(\mathbf{Y}_1,\mathbf{Y}_2,\cdots,\mathbf{Y}_p)$. Also, $\mathbf{x}=(\mathbf{s}_1^T,\mathbf{s}_2^T,\cdots,\mathbf{s}_p^T)^T$, and  $\{\mathbf{\tilde a}\}=(\mathbf{\tilde a}_1^T,\mathbf{\tilde a}_2^T,\cdots,\mathbf{\tilde a}_p^T)^T$ where $\mathbf{\tilde{a}}_i$ denotes an error of the parameter estimate
such that $\mathbf{\tilde{a}}_i=\mathbf{\hat{a}}_i-\mathbf{a}_i$. Note
that $\mathbf{a}_i$ is a constant vector of the true parameter values for the $i$-th robot,
resulting in $\mathbf{\dot{\tilde a}}_i=\mathbf{\dot{\hat a}}_i$. If each robot is identical, $\mathbf{a}_i=\mathbf{a}$ for $1\leq i \leq p$.
\begin{theorem}
\emph{The adaptive synchronization law in
(\ref{adaptive_sync_law}) globally asymptotically synchronizes the states of
multiple dynamics in the presence of parametric model
uncertainties if the condition (\ref{semi_contraction_condition}) holds.}
\end{theorem}
\begin{proof}
Similar to Section \ref{mainproof}, applying the spectral
transformation, using the augmented $
\mathbf{V}_a=\text{diag}\left(
\mathbf{V} ,\mathbf{I}_{pn}\right)$ and $ \mathbf{V}^T[\mathbf{L}^{p}_{\mathbf{K}_1,-\mathbf{K}_2}]\mathbf{V}=[\mathbf{D}]$, to (\ref{adaptive_block}) leads to the following virtual system of $(\mathbf{y}_1^T,\mathbf{y}_2^T)^T$~\cite{Jouffroy,Ref:contraction2}
\begin{align}\label{adaptive_block_virtual}
\begin{bmatrix}\mathbf{V}^T[\mathbf{M}]\mathbf{V} & \mathbf{0}\\ \mathbf{0} & [\mathbf{\Gamma^{-1}}] \end{bmatrix}&\begin{pmatrix}\mathbf{\dot y}_1\\ \mathbf{\dot y}_2\end{pmatrix}
      +\begin{bmatrix} \mathbf{V}^T[\mathbf{C}]\mathbf{V}&\mathbf{0} \\ \mathbf{0}&\mathbf{0} \end{bmatrix}\begin{pmatrix}\mathbf{ y}_1\\ \mathbf{ y}_2\end{pmatrix}\\&+\begin{bmatrix}[\mathbf{D}] & -\mathbf{V}^T[\mathbf{Y}] \\ [\mathbf{Y}]^T\mathbf{V}& \mathbf{0}\end{bmatrix}\begin{pmatrix}\mathbf{ y}_1\\ \mathbf{ y}_2\end{pmatrix}=\mathbf{0}\nonumber
\end{align}
The virtual system has two particular solutions:\\$
\mathbf{y}_1=\mathbf{V}^T\mathbf{x},\ \mathbf{y}_2=\{\mathbf{{\tilde a}}\}, \ \ \text{and} \ \ \ \mathbf{y}_1=\mathbf{0},\ \mathbf{y}_2=\mathbf{0}$\\
The virtual length analysis indicates that (\ref{adaptive_block_virtual}) is semi-contracting by the negative semi-definite Jacobian with $[\mathbf{D}]>0$:
\begin{equation}
\frac{dV}{dt}=-2\begin{pmatrix}\delta \mathbf{ y}_1\\ \delta \mathbf{ y}_2\end{pmatrix}^T\begin{bmatrix}[\mathbf{D}] & \mathbf{0} \\ \mathbf{0}& \mathbf{0}\end{bmatrix}\begin{pmatrix}\delta \mathbf{ y}_1\\ \delta \mathbf{ y}_2\end{pmatrix}\end{equation}
Using Barbalat's lemma (see Section \ref{sec:sync_no_tracking}), it
is straightforward to show that $\delta \mathbf{y}_1$
tends asymptotically to zero from any initial condition. Also, $[\mathbf{D}]$ can be decomposed to the tracking and synchronization gains, and the rest of the proof follows Theorems~\ref{Thm:tracking_sync} and \ref{Thm:tracking_sync_unstable}. Consequently,
the adaptive synchronization law in
(\ref{adaptive_sync_law}) synchronizes the states of
multiple dynamics in the presence of parametric model
uncertainties.\end{proof}

While the synchronization of the estimates of the physical parameters ($\delta \mathbf{y}_2\rightarrow \mathbf{0}$) is not automatically guaranteed due to the semi-contracting stability, the additional condition of the persistency of excitation~\cite{Ref:Slotine} leads to the convergence of $\mathbf{{\tilde a}}$ to zero. Results of the simulation are presented in Section~\ref{sec:adaptive_sim}.
\subsection{Examples}\label{two_robot_example}
For the case of a two-robot network, we can easily verify\begin{align}
&\mathbf{z}=\mathbf{V}^T\mathbf{x}=\begin{pmatrix}[\mathbf{1}]^T\\ \mathbf{V}_{sync}^T\end{pmatrix}\mathbf{x}=\begin{bmatrix} \frac{1}{\sqrt{2}}\mathbf{I} &
\frac{1}{\sqrt{2}}\mathbf{I}\\ \frac{1}{\sqrt{2}}\mathbf{I} &
-\frac{1}{\sqrt{2}}\mathbf{I}\end{bmatrix}\begin{pmatrix}\mathbf{s}_1\\\mathbf{s}_2\end{pmatrix},\nonumber \\ &\mathbf{D}_1=\mathbf{K}_1-\mathbf{K}_2,\ \ \mathbf{D}_2=\mathbf{K}_1+\mathbf{K}_2\end{align}

If $\mathbf{K}_1+\mathbf{K}_2>\mathbf{K}_1-\mathbf{K}_2>0$, the rate
of the virtual length $V$ in (\ref{virtual_length_good}) is uniformly
negative definite:
\begin{align}
\dot{V}=
-2\begin{pmatrix}\delta\mathbf{y}_t\\
\delta\mathbf{y}_s\end{pmatrix}^T\begin{bmatrix}
\mathbf{K}_1-\mathbf{K}_2
& \mathbf{0} \\ \mathbf{0} & \mathbf{K}_1+\mathbf{K}_2\end{bmatrix} \begin{pmatrix}\delta\mathbf{y}_t\\
\delta\mathbf{y}_s\end{pmatrix}       <0.\nonumber
\end{align}

\noindent Consequently, the combined virtual system in (\ref{virtual_chap5})
is contracting. As a result,
$\mathbf{s}_1+\mathbf{s}_2\rightarrow 0$ and
$\mathbf{s}_1\rightarrow\mathbf{s}_2$ exponentially.
It is straightforward to show that $\mathbf{s}_1\rightarrow\mathbf{s}_2$ also
hierarchically makes $\mathbf{q}_1$ tend to $\mathbf{q}_2$
exponentially as in (\ref{q_contraction}).

%

The results hold for arbitrarily large
networks as well. For example, a network of three robots has the following
$\mathbf{V}$, whose columns are orthonormal eigenvectors of
$[\mathbf{L}^{p=3}_{\mathbf{K}_1,-\mathbf{K}_2}]$:\begin{equation*}
\mathbf{V}=\begin{bmatrix}[\mathbf{1}] & \mathbf{V}_{sync}\end{bmatrix}=\left[\begin{smallmatrix}\frac{1}{\sqrt{3}}\mathbf{I} &
-\frac{2}{\sqrt{6}}\mathbf{I} & \mathbf{0}\\
\frac{1}{\sqrt{3}}\mathbf{I} & \frac{1}{\sqrt{6}}\mathbf{I} &
-\frac{1}{\sqrt{2}}\mathbf{I}\\ \frac{1}{\sqrt{3}}\mathbf{I} &
\frac{1}{\sqrt{6}}\mathbf{I} &
\frac{1}{\sqrt{2}}\mathbf{I}\end{smallmatrix}\right]\end{equation*}The block diagonal matrix $[\mathbf{D}]$ is also computed as $\text{diag}\left(\mathbf{K}_1-2\mathbf{K}_2 ,
\mathbf{K}_1+\mathbf{K}_2 ,
\mathbf{K}_1+\mathbf{K}_2\right)$.

For a four-robot network $p=4$ (see Figure~\ref{robot_diagram}(a)), $[\mathbf{D}]=\text{diag}\left(\mathbf{K}_1-2\mathbf{K}_2 ,
\mathbf{K}_1+2\mathbf{K}_2 ,\mathbf{K}_1,\mathbf{K}_1\right)$, and \begin{equation*}\mathbf{V}=\frac{1}{2}\left[\begin{smallmatrix}\mathbf{I}&\mathbf{I}&\mathbf{0}&-\sqrt{2}\mathbf{I}\\
\mathbf{I}&-\mathbf{I}&-\sqrt{2}\mathbf{I}&\mathbf{0}\\ \mathbf{I}&\mathbf{I}&\mathbf{0}&\sqrt{2}\mathbf{I}\\ \mathbf{I}&-\mathbf{I}&\sqrt{2}\mathbf{I}&\mathbf{0}\end{smallmatrix}\right]
\end{equation*}By inspecting the associated eigenvectors, we can notice that $\mathbf{K}_1$ represents the synchronization gain associated with the synchronization of diagonal members ($\mathbf{q}_1=\mathbf{q}_3,\mathbf{q}_2=\mathbf{q}_4$), while $\mathbf{K}_1+2\mathbf{K}_2$ represents the synchronization gain of direct couplings (e.g., $\mathbf{q}_1=\mathbf{q}_2,\mathbf{q}_3=\mathbf{q}_4$). This is a percolation effect discussed in \cite{Ref:contraction_sync}. The percolation effect can be exploited in order to prove the synchronization to the large unknown invariant set (e.g., $\mathbf{V}_{sync}^T\mathbf{x}=0$ from a large network) by verifying the synchronization to a known, not necessarily invariant, subset of the global flow-invariant set~\cite{Ref:Gerard_Slotine}.
\section{Concurrent Synchronization of Heterogeneous Groups on Unbalanced Graphs}\label{sec:concurrent}
We further generalize the proposed synchronization framework in the context of the concurrent synchronization of multiple heterogeneous networks and leader-follower networks, which permit construction of complex network structures. 
\subsection{Inline Configuration and Uni-Directional Digraph}\label{sec:inline}
We can also consider the inline configuration shown in Fig.~\ref{robot_diagram}(f), when maintaining a ring structure is not feasible possibly due to communication problems.
\begin{corollary}\label{CorollaryInline}
\emph{All the previous theorems hold for network structures on inline configuration, if the first and last robot adjust the proposed control law in (\ref{tracking_controllerCH5}) as
\begin{equation}\label{tracking_controllerInline}
\mathbf{\tau}_i={\mathbf{M}} (\mathbf{q}_i){\mathbf{\ddot{q}}_{i,r}}+
       \mathbf{C}(\mathbf{q}_i,\mathbf{\dot{q}}_i){\mathbf{\dot{q}}_{i,r}}+\mathbf{g}(\mathbf{q}_i)-(\mathbf{K}_1-\mathbf{K}_2)\mathbf{s}_i+\mathbf{K}_2\mathbf{s}_{j}
\end{equation}
where $i=1$ or $i=p$, and $j$ is the index for its sole neighbor.}
\end{corollary}
\begin{proof}
We can show that the modified Laplacian $[\mathbf{L}^{p}_{\mathbf{K}_1,-\mathbf{K}_2}]$ constructed from both (\ref{tracking_controllerCH5}) and (\ref{tracking_controllerInline}) is symmetric and still has $[\mathbf{1}]$ as its eigenvector associated its eigenvalue $\mathbf{D}_1=\mathbf{K}_1-2\mathbf{K}_2$. Then the proofs of the theorems easily follow.
\end{proof}

Let us now consider the generalized control law in (\ref{tracking_controllerCH5_general}), which permits regular \emph{digraphs}, as shown in Figures~\ref{robot_diagram}(c) and (d). The modified Laplacian matrix $[\mathbf{L}^{p}_{\mathbf{K}_1,-\mathbf{K}_2}]$ defined from (\ref{tracking_controllerCH5_general}) might have two or more than three nonzero elements in each row. Further, it is no longer symmetric.
\begin{corollary}
\emph{All the previous theorems are valid for regular graphs with uni-directional couplings or both uni-directional and bi-directional couplings.}
\end{corollary}
\begin{proof}
It is straightforward to verify that the symmetric matrix $[\mathbf{L}^{p}_{\mathbf{K}_1,-\mathbf{K}_2}]+[\mathbf{L}^{p}_{\mathbf{K}_1,-\mathbf{K}_2}]^T$ has the same number of nonzero elements at each row (i.e., regular). Hence, we can conclude that all the previous proofs still hold. For example, the synchronization occurs when
$\mathbf{V}^T_{sync}\left([\mathbf{L}^{p}_{\mathbf{K}_1,-\mathbf{K}_2}]+[\mathbf{L}^{p}_{\mathbf{K}_1,-\mathbf{K}_2}]^T\right)\mathbf{V}_{sync}>0
$ while $[\mathbf{1}]^T\left([\mathbf{L}^{p}_{\mathbf{K}_1,-\mathbf{K}_2}]+[\mathbf{L}^{p}_{\mathbf{K}_1,-\mathbf{K}_2}]^T\right)[\mathbf{1}]
$ determines the stability and convergence rate of the trajectory tracking.\end{proof}
\subsection{Synchronization of Heterogeneous Robots}\label{sec:nonid_sync}
Consider the proposed control law for a network comprised of heterogeneous dynamics as follows:
\begin{align}
&  \mathbf{\tau}_i={\mathbf{M}}_i (\mathbf{q}_i){\mathbf{\ddot{q}}_{i,r}}+
       \mathbf{C}_i(\mathbf{q}_i,\mathbf{\dot{q}}_i){\mathbf{\dot{q}}_{i,r}}+\mathbf{g}_i(\mathbf{q}_i)
      \\& -\mathbf{K}_1(\mathbf{\dot{q}}_i-\mathbf{\dot{q}}_{i,r})+\mathbf{K}_2(\mathbf{\dot{q}}_{i-1}-\mathbf{\dot{q}}_{i-1,r})+\mathbf{K}_2(\mathbf{\dot{q}}_{i+1}-\mathbf{\dot{q}}_{i+1,r})\nonumber
\end{align}
where ${\mathbf{M}}_i$, $\mathbf{C}_i$, and
$\mathbf{g}_i(\mathbf{q}_i)$ represent the $i$-th robot dynamics,
which can be different from robot to robot. Each robot has the same
number of configuration variables ($\mathbf{q}_i\in\mathbb{R}^n$).
\begin{corollary}
\emph{The proposed tracking and synchronization control law in
(\ref{tracking_controllerCH5}) can easily be applied to a network
consisting of heterogeneous robots in (\ref{NL_single_compactCH5}) if
the stable tracking condition in Theorem \ref{Thm:tracking_sync} is
true.}\end{corollary}
\begin{proof}
The $\mathbf{M}$ and $\mathbf{C}$ matrix
notations used in the previous sections can be interpreted as
$\mathbf{M}(\mathbf{q}_1)\rightarrow\mathbf{M}_1(\mathbf{q}_1)$ and
$\mathbf{M}(\mathbf{q}_2)\rightarrow\mathbf{M}_2(\mathbf{q}_2)$ with
$\mathbf{M}_1(\cdot)\neq\mathbf{M}_2(\cdot)$ (the same for the
$\mathbf{C}$ matrices). Hence, the assumption of non-identical
dynamics does not alter the proof of exponential tracking in Theorem
\ref{Thm:tracking_sync_conv} and stable synchronization in Theorem~\ref{Thm:tracking_sync}. However, the synchronization with indifferent tracking (Theorem \ref{Thm:tracking_sync_unstable}) is no longer true
for non-identical robots since $\mathbf{V}^T_{sync}\mathbf{x}=0$ is not a flow-invariant manifold, and
$\mathbf{q}_1=\mathbf{q}_2$ does not cancel the off-diagonal terms
$1/2\left[{\mathbf{M}_1} (\mathbf{q}_1)-{\mathbf{M}_2}
(\mathbf{q}_2)\right]$ in the metric matrix $\mathbf{V}^T[\mathbf{M}]\mathbf{V}$.\end{proof}
\subsection{Concurrent Synchronization of Heterogeneous Groups}
We present a new method of achieving the concurrent synchronization of multiple heterogeneous networks in this section. Such a method can be used to apply the main results in the previous sections to more complex networks, as shown in Figure~\ref{robot_diagram}(e). The networks in the figure are \emph{neither regular nor balanced} due to the reference input couplings.  In \cite{Ref:contraction_sync}, \emph{concurrent synchronization} is
defined as a regime where the ensemble of dynamical elements is
divided into multiple groups of fully synchronized elements, but
elements from different groups are not necessarily synchronized and
can exhibit entirely different dynamics.
\begin{figure}
 \begin{center}
  \includegraphics[width=3.5in]{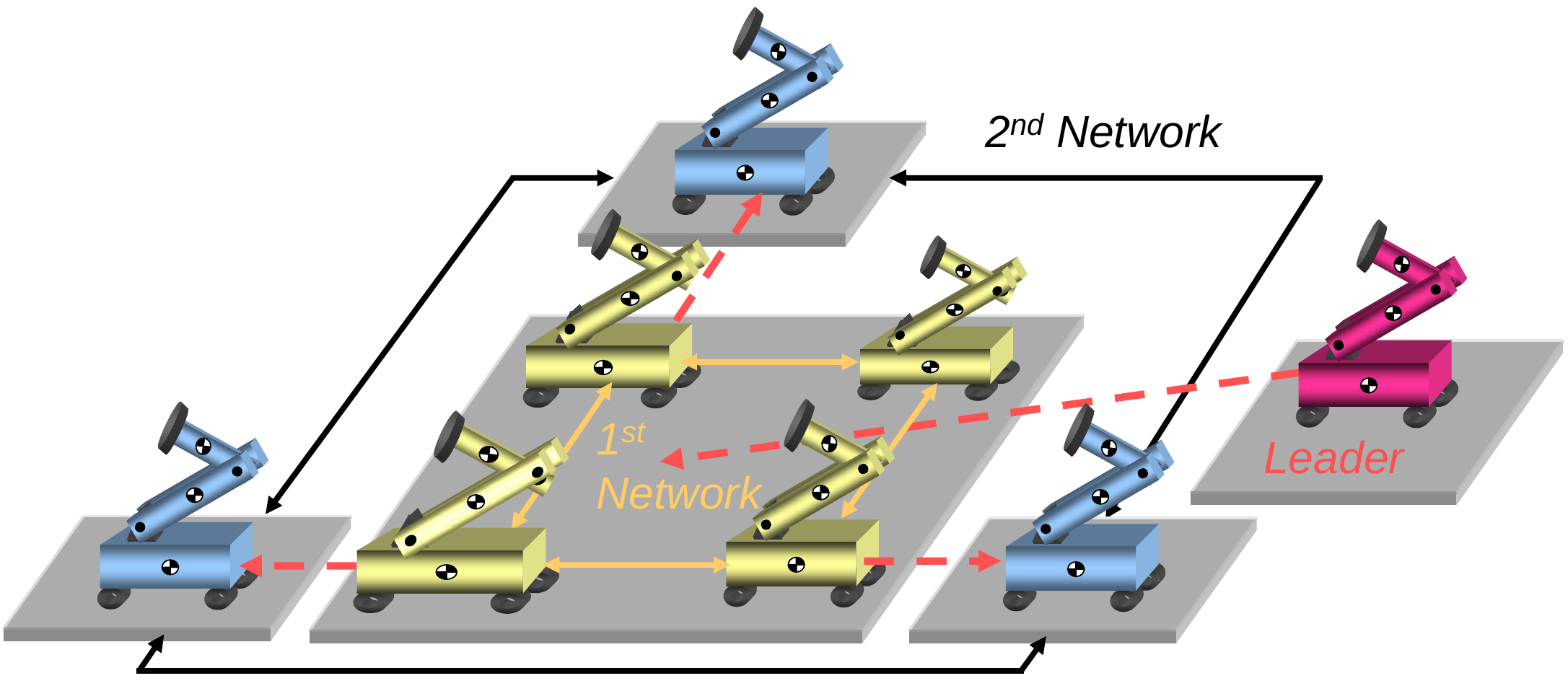}\end{center}
  \caption{Concurrent synchronization between two different groups. The desired trajectory inputs are denoted by the dashed-lines whereas the solid lines indicate mutual diffusive couplings. The independent leader sends the same desired trajectory input to the first network group. If we view the dashed-lines as edges of the graph, the network is on an \emph{unbalanced} graph.}\label{concurrent_sync}
\end{figure}

As discussed in the previous
sections, we pay particular attention to the fact that there exist
two different time scales of the proposed synchronization tracking
control law. In particular, one rate is
associated with the trajectory tracking
($\mathbf{D}_1$), while the other represents the
convergence rate of the synchronization
($\mathbf{D}_2$). This implies that there are
two different inputs to the system, namely, the common reference trajectory $\mathbf{q}_d(t)$ and the diffusive couplings with the
adjacent members. Accordingly, we exploit a desired trajectory $\mathbf{q}_d(t)$ to
create multiple combinations of different dynamics groups.

For
instance, Figure \ref{concurrent_sync} represents the concurrent
synchronization of two different dynamic networks. The first
network, consisting of four heterogeneous robots, has the diffusive
coupling structure proposed by the tracking control law in
(\ref{tracking_controllerCH5}). The independent leader sends a
desired trajectory command $\mathbf{q}_d(t)$ to each member of the first network. Note that the dimension of the leader dynamics need not be equal to that of the first network, since a signal from the leader can be pre-conditioned. With
an appropriate selection of gains, each dynamics in the first
network synchronize while exponentially following the leader.
Therefore, the proposed scheme can be interpreted in the context of
the leader-follower problem
(\cite{Jadbabaie,Ref:Leonard_leader,Ref:contraction_leader}).

The second network consists of three heterogeneous dynamics, also
different from those of the first group. Each element receives a
different desired trajectory input from the adjacent element of the
first network. Again, the desired trajectory input for the second group can be pre-conditioned from the first network, resulting in
\begin{equation}
\mathbf{q}_{2d}(t)=\mathbf{f}(\mathbf{q}_1,t)
\end{equation}
where $\mathbf{q}_1$ is the state vector of the member in the first network that sends in the desired command and $\mathbf{f}$ is a differentiable nonlinear function of $\mathbf{q}_1$. As a result, the dynamics of the two networks can be entirely different (i.e., $\mathbf{q}_1\in\mathbb{R}^m$, $\mathbf{q}_2\in\mathbb{R}^n$, and $m\neq n$). Note that $\mathbf{\ddot q}_{2d}$, needed for the proposed control law in (\ref{tracking_controllerCH5_general}), can be obtained from a nonlinear observer~\cite{Ref:Slotine}. The concurrent synchronization between the first and second networks can be proven as follows.
\begin{corollary}\label{Corollary_con}
\emph{Consider the network structures that are unbalanced (i.e., the in-degrees are not equal to the out-degrees) due to the directional reference input connections among the regular balanced graphs, shown as the dashed lines in Figures~\ref{robot_diagram}(e) and \ref{concurrent_sync}. The robots on such graphs globally exponentially synchronize if the individual groups synchronize from Theorems~\ref{Thm:tracking_sync_conv} and \ref{Thm:tracking_sync}.}
\end{corollary}
\begin{proof}
For the example shown in Figure~\ref{concurrent_sync}, once the first network synchronizes, the
second network also ends up receiving the same desired trajectory to
follow, while they interact to synchronize exponentially fast. Accordingly, we can achieve concurrent synchronization between two different network groups. The proof of Theorem~\ref{Thm:tracking_sync} holds until (\ref{q_contraction}), where now the desired trajectory inputs are different for each robot. Then, we can conclude that the robots in the second network synchronize, if the desired reference inputs $\mathbf{q}_{2d,i}$, $1 \leq i \leq p$, sent from the first network, synchronize:
\begin{align}
&\mathbf{V}_{sync}^T\{\mathbf{\dot q}_2\}+\left(\mathbf{V}_{sync}^T[\mathbf{\Lambda}]\mathbf{V}_{sync}\right)\mathbf{V}_{sync}^T\{\mathbf{q}_2\}\\ \nonumber
&\ \  -\left(\mathbf{V}_{sync}^T[\mathbf{r}_{1}^T,\cdots,\mathbf{r}_{p}^T]^T\rightarrow \mathbf{0}\right)=\mathbf{V}_{sync}^T\mathbf{x}_2\rightarrow \mathbf{0}
\end{align}
where $\mathbf{r}_i=\mathbf{\dot q}_{2d,i}+\mathbf{\Lambda}\mathbf{q}_{2d,i}$, and $\mathbf{x}_2$ is the vector of the composite variables ($\mathbf{s}_i$) of the second network. If the first network synchronizes, the reference trajectory term vanishes ($\mathbf{V}_{sync}^T[\mathbf{r}_{1}^T,\cdots,\mathbf{r}_{p}^T]^T\rightarrow \mathbf{0}$), thereby ensuring the synchronization of the second network (i.e., $\mathbf{V}_{sync}^T\{\mathbf{q}_2\}\rightarrow \mathbf{0}$).
This can be extended to arbitrarily large
groups of synchronized dynamics by appropriately assigning the
desired trajectory inputs and the diffusive couplings.
\end{proof}
\begin{remark}
The previous example in Figure~\ref{concurrent_sync} and Corollary~\ref{Corollary_con} can be interpreted in terms of feedback hierarchies (see Theorem~\ref{Thm:hierc}). The first network in the figure provides feedforward commands to the second network, but does not receive any command from the second network. Note that the dynamics at each level can be very different, in order to construct a dynamic concurrent combination of heterogeneous networks. Since the feedforward inputs can be appropriately scaled and conditioned, the dimensions between hierarchical layers can also be very different. For example, the dynamics higher in the hierarchy need not be oscillators, but could be systems with multiple equilibria (e.g., $\mathbf{\dot x} = -\nabla V$, with local minima in $V$~\cite{Ref:contraction2}), with synchronization corresponding to convergence toward a common minimum. Concurrent synchronization discussed in this section can be exploited to construct a large complex network consisting of heterogeneous dynamics such as robots, ground vehicles, and unmanned aerial vehicles.\end{remark}
\section{Extensions}\label{sec:ext_sim}
Let us further extend the proposed control law.
\subsection{Synchronization of Robots Using Linear PD Control}\label{SEC:PD_Sync_Big}
One may consider the following Proportional and Derivative (PD)
coupling control law for two identical robots from
(\ref{NL_single_compactCH5}) with $p=2$,
\begin{equation}\label{tracking_controllerpd}
\begin{split}
  \mathbf{\tau}_1=-\mathbf{K}_1(\mathbf{\dot{q}}_1+\mathbf{\Lambda}\mathbf{\tilde q}_1)+\mathbf{K}_2(\mathbf{\dot{q}}_2+\mathbf{\Lambda}\mathbf{\tilde q}_2)\\
   \mathbf{\tau}_2=-\mathbf{K}_1(\mathbf{\dot{q}}_2+\mathbf{\Lambda}\mathbf{\tilde q}_2)+\mathbf{K}_2(\mathbf{\dot{q}}_1+\mathbf{\Lambda}\mathbf{\tilde q}_1)
\end{split}
\end{equation}
where $\mathbf{\tilde q}_i=\mathbf{q}_i-\mathbf{q}_{d}$ and the
bounded reference position $\mathbf{q}_d$ has zero velocity such
that $\mathbf{\dot{\tilde q}}_i=\mathbf{\dot q}_i$.

\noindent For simplicity, the gravity term in (\ref{NL_single_compactCH5}) is
assumed to be zero or canceled by a feed-forward control law. Then,
the closed-loop dynamics satisfy
\begin{equation}\label{closed_single}
\begin{split}
  {\mathbf{M}} (\mathbf{q}_1){\mathbf{\ddot{q}}_{1}}+
       \mathbf{C}(\mathbf{q}_1,\mathbf{\dot{q}}_1){\mathbf{\dot{q}}_{1}}+\mathbf{K}\mathbf{\dot{q}}_1+\mathbf{K}\mathbf{\Lambda}\mathbf{\tilde{q}}_1=\mathbf{u}(t)\\
{\mathbf{M}} (\mathbf{q}_2){\mathbf{\ddot{q}}_{2}}+
       \mathbf{C}(\mathbf{q}_2,\mathbf{\dot{q}}_2){\mathbf{\dot{q}}_{2}}+\mathbf{K}\mathbf{\dot{q}}_2+\mathbf{K}\mathbf{\Lambda}\mathbf{\tilde{q}}_2=\mathbf{u}(t)
\end{split}
\end{equation}
where $\mathbf{K}=\mathbf{K}_1+\mathbf{K}_2$ and
$\mathbf{u}(t)=\mathbf{K}_2(\mathbf{\dot{q}}_1+\mathbf{\dot{q}}_2)+\mathbf{K}_2\mathbf{\Lambda}(\mathbf{\tilde{q}}_1+\mathbf{\tilde{q}}_2)$.

Similar to Section~\ref{mainproof}, we can perform a spectral
decomposition:
\begin{align}
&\left(\mathbf{V}^T[\mathbf{M}]\mathbf{V}\right)\mathbf{V}^T\mathbf{\ddot
x}+       \left(\mathbf{V}^T[\mathbf{C}]\mathbf{V}\right)\mathbf{V}^T\mathbf{\dot
x} \\&+\left(\mathbf{V}^T[\mathbf{L}^{p}_{\mathbf{K}_1,-\mathbf{K}_2}]\mathbf{V}\right)\mathbf{V}^T\mathbf{\dot
x}+\left(\mathbf{V}^T[\mathbf{L}^{p}_{\mathbf{K}_1\mathbf{\Lambda},-\mathbf{K}_2\mathbf{\Lambda}}]\mathbf{V}\right)\mathbf{V}^T\mathbf{\tilde
x}=\mathbf{0}\nonumber
\end{align}

Using the following Lyapunov function, it is straightforward to show
that this PD coupling control law drives the system to the desired
rest state $\mathbf{q}_d$ globally and asymptotically while tending
to the synchronized flow-invariant manifold:
\begin{align}
V=&\frac{1}{2}\begin{pmatrix}\mathbf{\dot q}_p\\
\mathbf{\dot
q}_m\end{pmatrix}^T\begin{bmatrix}\frac{{\mathbf{M}}_1+{\mathbf{M}}_2}{2}
&
\frac{{\mathbf{M}}_1-{\mathbf{M}}_2}{2}\\
\frac{{\mathbf{M}}_1-{\mathbf{M}}_2}{2}&\frac{{\mathbf{M}}_1+{\mathbf{M}}_2}{2}\end{bmatrix}\begin{pmatrix}\mathbf{\dot q}_p\\
\mathbf{\dot q}_m\end{pmatrix}\\&+\frac{1}{2}\begin{pmatrix}\mathbf{\tilde q}_p\\
\mathbf{q}_m\end{pmatrix}^T\begin{bmatrix}(\mathbf{K}_1-\mathbf{K}_2)\mathbf{\Lambda} & 0\\0 & (\mathbf{K}_1+\mathbf{K}_2)\mathbf{\Lambda}\end{bmatrix}\begin{pmatrix}\mathbf{\tilde q}_p\\
\mathbf{q}_m\end{pmatrix}\nonumber
\end{align}
where ${\mathbf{M}}_1={\mathbf{M}}
(\mathbf{q}_1),\ \ {\mathbf{M}}_2={\mathbf{M}} (\mathbf{q}_2),\ \
\mathbf{V}^T\mathbf{ x}=\begin{pmatrix}\mathbf{q}_p
\\
\mathbf{q}_m\end{pmatrix}=\frac{1}{\sqrt{2}}\begin{pmatrix}\mathbf{q}_1+\mathbf{q}_2
\\ \mathbf{q}_1-\mathbf{q}_2\end{pmatrix},\ \ \ \mathbf{\tilde
q}_p=\frac{1}{\sqrt{2}}(\mathbf{\tilde q}_1+\mathbf{\tilde
q}_2).$

The rate of $V$ can be computed as
\begin{equation}
\begin{split}
\frac{dV}{dt}=-\begin{pmatrix}\mathbf{\dot q}_p\\
\mathbf{\dot q}_m\end{pmatrix}^T\begin{bmatrix}\mathbf{K}_1-\mathbf{K}_2 & 0\\0 & \mathbf{K}_1+\mathbf{K}_2\end{bmatrix}\begin{pmatrix}\mathbf{\dot q}_p\\
\mathbf{\dot q}_m\end{pmatrix}\leq 0
\end{split}
\end{equation}
which implies that $\dot V$ is negative semi-definite with
$\mathbf{K}_1+\mathbf{K}_2>0$ and $\mathbf{K}_1-\mathbf{K}_2>0$.

By invoking LaSalle's invariant set theorem~\cite{Ref:Slotine}, we
can conclude that $\mathbf{\dot{{q}}}_p$, $\mathbf{\tilde q}_p$,
$\mathbf{\dot q}_m$, and $\mathbf{q}_m$ tend to zero with global and asymptotic convergence. This implies that
$\mathbf{q}_1$ and $\mathbf{q}_2$ will follow $\mathbf{q}_d$ while
$\mathbf{q}_1$ and $\mathbf{q}_2$ synchronize asymptotically from any initial condition. This can be extended to arbitrarily large networks as shown in Figure~\ref{robot_diagram}.

Note that if $\mathbf{\Lambda}=0$, the PD coupling control law in
(\ref{tracking_controllerpd}) reduces to velocity coupling
\begin{equation}\label{tracking_controllerpd_simple}
\begin{split}
  \mathbf{\tau}_1=-\mathbf{K}_1\mathbf{\dot{q}}_1+\mathbf{K}_2\mathbf{\dot{q}}_2, \ \
   \mathbf{\tau}_2=-\mathbf{K}_1\mathbf{\dot{q}}_2+\mathbf{K}_2\mathbf{\dot{q}}_1
\end{split}
\end{equation}

\noindent This velocity coupling control can also be derived from the
exponential tracking control law in (\ref{tracking_controllerCH5})
by setting $\mathbf{q}_d=0$, $\mathbf{\dot q}_d=0$, and
$\mathbf{\Lambda}=0$. Therefore, the proof of the linear PD synchronization with $\mathbf{\Lambda}=0$ is the same as Section~\ref{mainproof} whereas the convergence rate is now exponential
compared with the asymptotic convergence of the PD control. On the
other hand, we can find that positions do not synchronize \emph{in
the absence of the gravity term} even though the velocities
synchronize exponentially fast.

\subsection{Synchronization with Limited Communication Bandwidth}
We now consider multiple dynamics with partially coupled joints (or
partially coupled variables). For example, we can assume that only
the lower joint is coupled in a two-robot system having two joint
variables with $\mathbf{q}=(x_1, x_2)^T$ for $(i=1,j=2)$ or $(i=2,j=1)$:
\begin{align}
\mathbf{\tau}_i={\mathbf{M}} (\mathbf{q}_i){\mathbf{\ddot{q}}_{i,r}}&+
       \mathbf{C}(\mathbf{q}_i,\mathbf{\dot{q}}_i){\mathbf{\dot{q}}_{i,r}}+\mathbf{g}(\mathbf{q}_i)
-\mathbf{K}_1\mathbf{s}_i\nonumber\\&+\mathbf{K}_2\begin{pmatrix}\dot{\tilde x}_1\\ 0\end{pmatrix}_{\mathbf{q}_{j}}+\mathbf{K}_2\mathbf{\Lambda}\begin{pmatrix}\widetilde{x}_1\\ 0\end{pmatrix}_{\mathbf{q}_{j}}
\end{align}

Nevertheless, Theorems \ref{Thm:tracking_sync_conv} and
\ref{Thm:tracking_sync} are true with diagonal matrices,
$\mathbf{K}_1$, $\mathbf{K}_2$ and $\mathbf{\Lambda}$, which can be
verified by writing the closed-loop system as in
(\ref{closed_single_tracking}):
\begin{align}\label{closed_single_tracking_partial}
  {\mathbf{M}} (\mathbf{q}_1)\mathbf{\dot s}_1+
       \mathbf{C}(\mathbf{q}_1,\mathbf{\dot{q}}_1)\mathbf{s}_1+(\mathbf{K}_1+\mathbf{K}_2\mathbf{P})\mathbf{s}_1=\mathbf{u}(t)\\ \nonumber
{\mathbf{M}} (\mathbf{q}_2)\mathbf{\dot s}_2+
       \mathbf{C}(\mathbf{q}_2,\mathbf{\dot{q}}_2)\mathbf{s}_2+(\mathbf{K}_1+\mathbf{K}_2\mathbf{P})\mathbf{s}_2=\mathbf{u}(t)\\ \nonumber
\mathbf{u}(t)=\mathbf{K}_2\mathbf{P}(\mathbf{s}_1+{\mathbf{s}_{2}}), \ \ \mathbf{P}=\text{diag}(1,0)
\end{align}
It is straightforward to prove that Theorems
\ref{Thm:tracking_sync_conv} and \ref{Thm:tracking_sync} still hold.
This is because $(\mathbf{K}_1+\mathbf{K}_2\mathbf{P})$ and $(\mathbf{K}_1-\mathbf{K}_2\mathbf{P})$ are still uniformly positive definite,
enabling exponential synchronization and exponential convergence to
the desired trajectory, respectively. Hence, we did not break any
assumption in the proof of Theorem \ref{Thm:tracking_sync}. This partial-state coupling also works for the semi-contracting case with (\ref{defsemicontractings}), as presented in Section \ref{sec:sync_no_tracking}. While the coupled variables synchronize to the weighted average of the initial conditions, the uncoupled variables synchronize to zero.
\subsection{Synchronization with Time-Delays}
The results may be extended to time-delayed couplings (see also ~\cite{Ref:Spong_Sync}). Consider for instance the two-robot dynamics (\ref{closed_single_tracking}), which becomes
\begin{align}\label{two_robot_time_new2}
&{\mathbf{M}} (\mathbf{q}_1)\mathbf{\dot s}_1+
       \mathbf{C}(\mathbf{q}_1,\mathbf{\dot{q}}_1)\mathbf{s}_1+(\mathbf{K}_1-\mathbf{K}_2)\mathbf{s}_1=\mathbf{\tau}_{21} \nonumber\\
&{\mathbf{M}} (\mathbf{q}_2)\mathbf{\dot s}_2+
       \mathbf{C}(\mathbf{q}_2,\mathbf{\dot{q}}_2)\mathbf{s}_2+(\mathbf{K}_1-\mathbf{K}_2)\mathbf{s}_2=\mathbf{\tau}_{12} \nonumber \\
&\mathbf{\tau}_{21}=\mathbf{K}_{2}(\mathbf{s}_2(t-T)-\mathbf{s}_1(t))\\ &\mathbf{\tau}_{12}=\mathbf{K}_{2}(\mathbf{s}_1(t-T)-\mathbf{s}_2(t))\nonumber
\end{align}
where the constant $T>0$ denotes the communication delay
between the two robots, and $\mathbf{s}(t-T)=\mathbf{\dot q}(t-T)+\mathbf{\Lambda} \mathbf{q}(t-T)-(\mathbf{\dot q}_d(t)+\mathbf{\Lambda} \mathbf{q}_d(t))$.
\begin{theorem}
\emph{The robots in (\ref{two_robot_time_new2}) synchronize globally asymptotically $\forall T>0$ under the assumptions of Theorem~\ref{Thm:tracking_sync}.}
\end{theorem}
\begin{proof}
Similar to~\cite{Ref:contraction5,Ref:contraction_timedelay}, consider
the  differential length
\begin{equation}
V=\delta \mathbf{z}^T \left(\mathbf{V}^T[\mathbf{M}]\mathbf{V}\right)\delta \mathbf{z}+\int^{t}_{t-T} \delta \mathbf{z}^T (\epsilon)\left[\begin{smallmatrix}\mathbf{K}_{2} & \mathbf{0}\\ \mathbf{0} & \mathbf{K}_{2}\end{smallmatrix}\right]\delta \mathbf{z} (\epsilon)d\epsilon\end{equation}
Using (\ref{two_robot_time_new2}) shows that $\dot V$ is negative semi-definite,
\begin{align}
\dot V=&-2\delta \mathbf{z}^T\left[\begin{smallmatrix}\mathbf{K}_{1}-\mathbf{K}_{2} & \mathbf{0}\\ \mathbf{0} & \mathbf{K}_{1}-\mathbf{K}_{2}\end{smallmatrix}\right]\delta\mathbf{z}\\&-\delta\mathbf{\tau}_{21}^T\mathbf{K}^{-1}_{2}\delta\mathbf{\tau}_{21}-\delta\mathbf{\tau}_{12}^T\mathbf{K}^{-1}_{2}\delta\mathbf{\tau}_{12}\nonumber
\end{align}
Note that $\mathbf{K}_1+\mathbf{K}_2>\mathbf{K}_1-\mathbf{K}_2>0$ implies $\mathbf{K}_2>0$.

Also, $\ddot V$ is bounded. By Barbalat's lemma, $\dot V$ tends to zero globally asymptotically, which in turn implies $\delta\mathbf{z}\rightarrow 0, \ \delta\mathbf{\tau}_{21}\rightarrow 0,\ \delta\mathbf{\tau}_{12}\rightarrow 0$. Hence, the solutions converge to a single trajectory asymptotically. Since $\mathbf{s}_1(t-T)=\mathbf{s}_2(t-T)$ is a particular solution of (\ref{two_robot_time_new2}), $\mathbf{q}_1(t)$ and $\mathbf{q}_2(t)$ globally asymptotically synchronize regardless of $T$.
\end{proof}
%
\subsection{A Perspective on Model Reduction by
Synchronization}\label{sec:sync_reduction}
Another benefit of exponential synchronization is its implication for model
reduction~\cite{Ref:tether}. Exponential synchronization of multiple
nonlinear dynamics allows us to reduce the dimensionality of the
stability analysis of a large network. As noted earlier, the
synchronization rate ($\mathbf{D}_2$) is faster than the
tracking rate ($\mathbf{D}_1$). Assuming that the
dynamics are synchronized, the augmented dynamics in
(\ref{tracking_form}) reduces to
\begin{equation}
\mathbf{M}(\mathbf{q})\mathbf{\dot s}+\mathbf{C}(\mathbf{q},\mathbf{\dot q}) \mathbf{s}+(\mathbf{K}_1-2\mathbf{K}_2)\mathbf{s}=\mathbf{0}
\end{equation}where $\mathbf{q}=\mathbf{q}_1=\cdots=\mathbf{q}_p$, and $\mathbf{D}_1=\mathbf{K}_1-2\mathbf{K}_2$ is replaced by $\mathbf{K}_1-\mathbf{K}_2$ for $p=2$.

\noindent This implies that once components of a network are shown to
synchronize, we can regard them as a single dynamics of reduced
dimension, which simplifies any additional stability or perturbation
analysis.
\section{Simulation Results}
\subsection{Tracking Synchronization of Four Robots}
\begin{figure}
 \begin{subfigmatrix}{1}
  \subfigure[Simulation result.]{\includegraphics[width=3.2in]{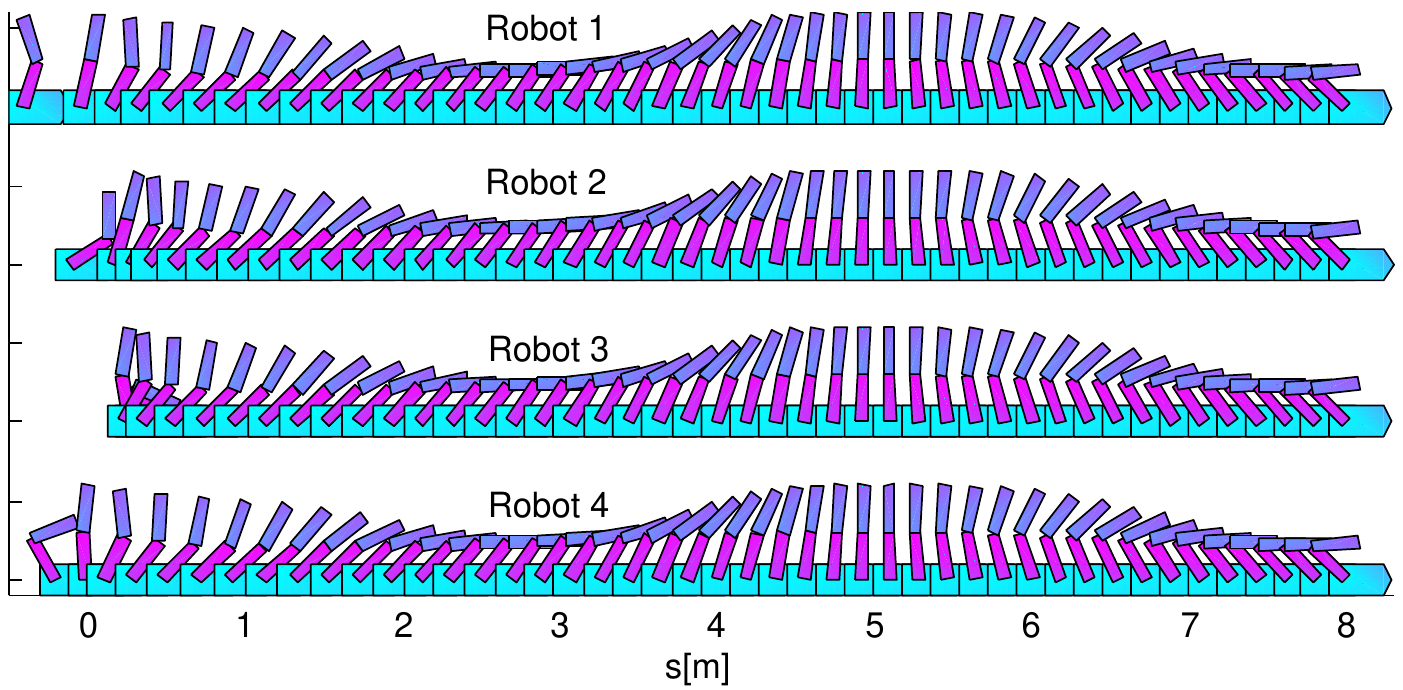}}
  \subfigure[Tracking error of each robot]{\includegraphics[width=3.1in]{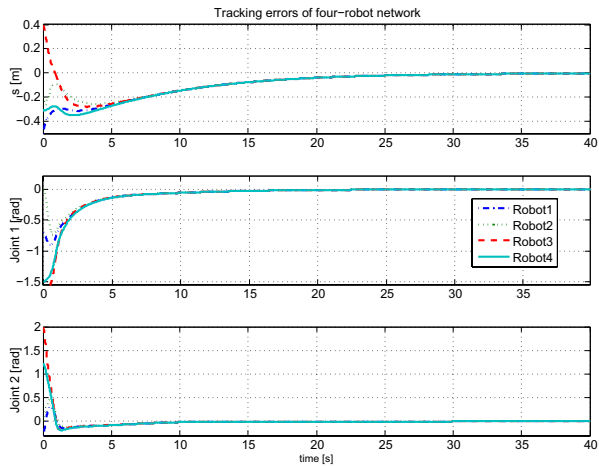}}
 \end{subfigmatrix}
 \caption{Synchronization of a four-robot network}\label{four_robot_sim}
\end{figure}

Even though the local coupling structure of (\ref{tracking_controllerCH5_general}) and (\ref{tracking_controllerCH5}) has been emphasized, the
difference from all-to-all coupling is not evident in a network
comprised of less than four members ($p\leq3$). To illustrate the
effectiveness of the proposed scheme for a robot network with $p\geq4$, a network of
four identical 3-DOF robots is considered here (see Fig~\ref{robot_diagram}(a)). The dynamics modeling of the 3-DOF robot is based upon the double
inverted pendulum robot on a cart (see Fig.~\ref{four_robot_sim}(a)), and each joint is assumed to be frictionless. The physical parameters of each robot are given in \cite{Ref:sjchung_PhD}.

The simulation result is presented in Figure \ref{four_robot_sim}.
The four identical robots, initially at some arbitrary initial conditions, are driven to synchronize as well as to
track the time-varying desired trajectory: $s_{d}(t)=0.2t,\ \theta_{1d}(t)=\cos(0.02\pi t),  \ \theta_{2d}(t)=\pi/4(1-\cos(0.08\pi t))$. For the control gains in the control law
(\ref{tracking_controllerCH5}), we used $\mathbf{K}_1=5\mathbf{I}$,
$\mathbf{K}_2=2\mathbf{I}$, and $\mathbf{\Lambda}=5\mathbf{I}$.
According to Theorems \ref{Thm:tracking_sync_conv} and
\ref{Thm:tracking_sync}, the tracking gain
$\mathbf{K}_1-2\mathbf{K}_2$ is smaller than the corresponding
synchronization gains $\mathbf{K}_1+2\mathbf{K}_2$ and $\mathbf{K}_1$. Figure
\ref{four_robot_sim}(b) shows the tracking errors of the four robots. We can see that the robots synchronize exponentially fast from arbitrary initial conditions, and this synchronization occurs faster than the exponential convergence of tracking errors. Such a result can be useful to rapidly achieve a collective motion of multiple robots in the presence of external disturbances. If $\mathbf{K}_1=2\mathbf{K}_2$ instead, we can easily show that the robots synchronize with global exponential convergence, although the tracking errors will not tend to zero.
\begin{figure}
 \begin{subfigmatrix}{1}
  \subfigure[Three networks connected by feedback hierarchies.]{\includegraphics[width=2.5in]{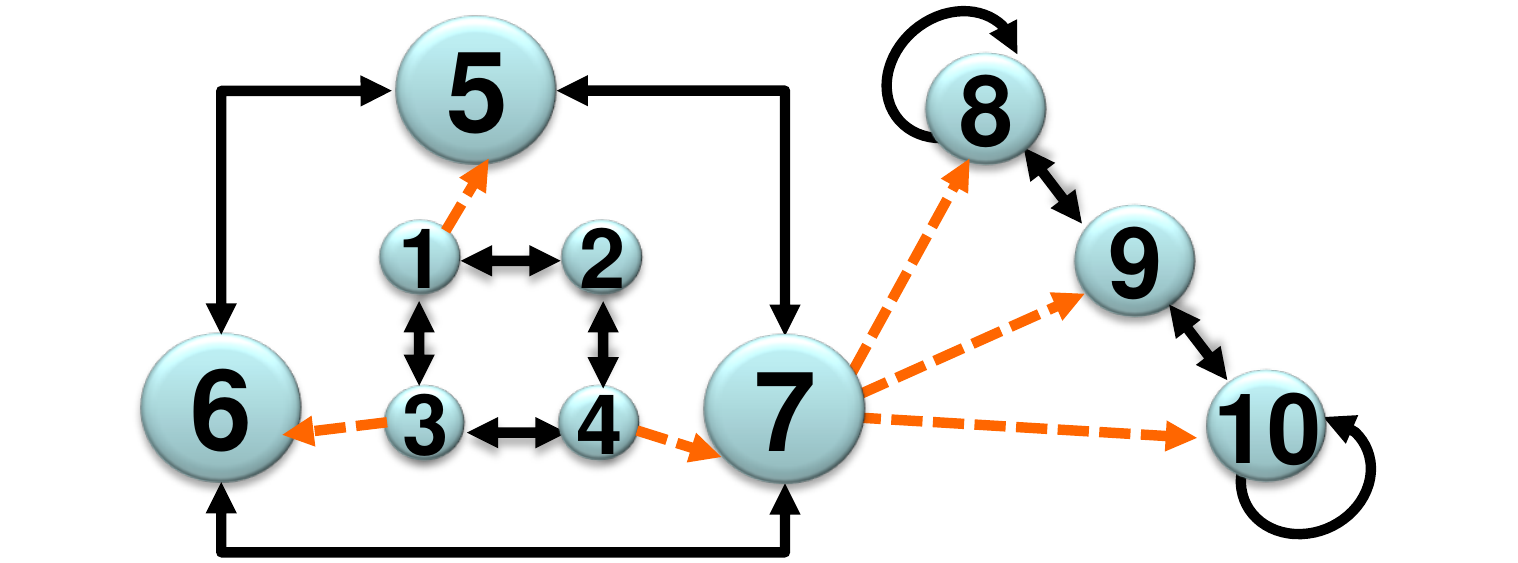}}
  \subfigure[Tracking error of each robot]{\includegraphics[width=3.5in]{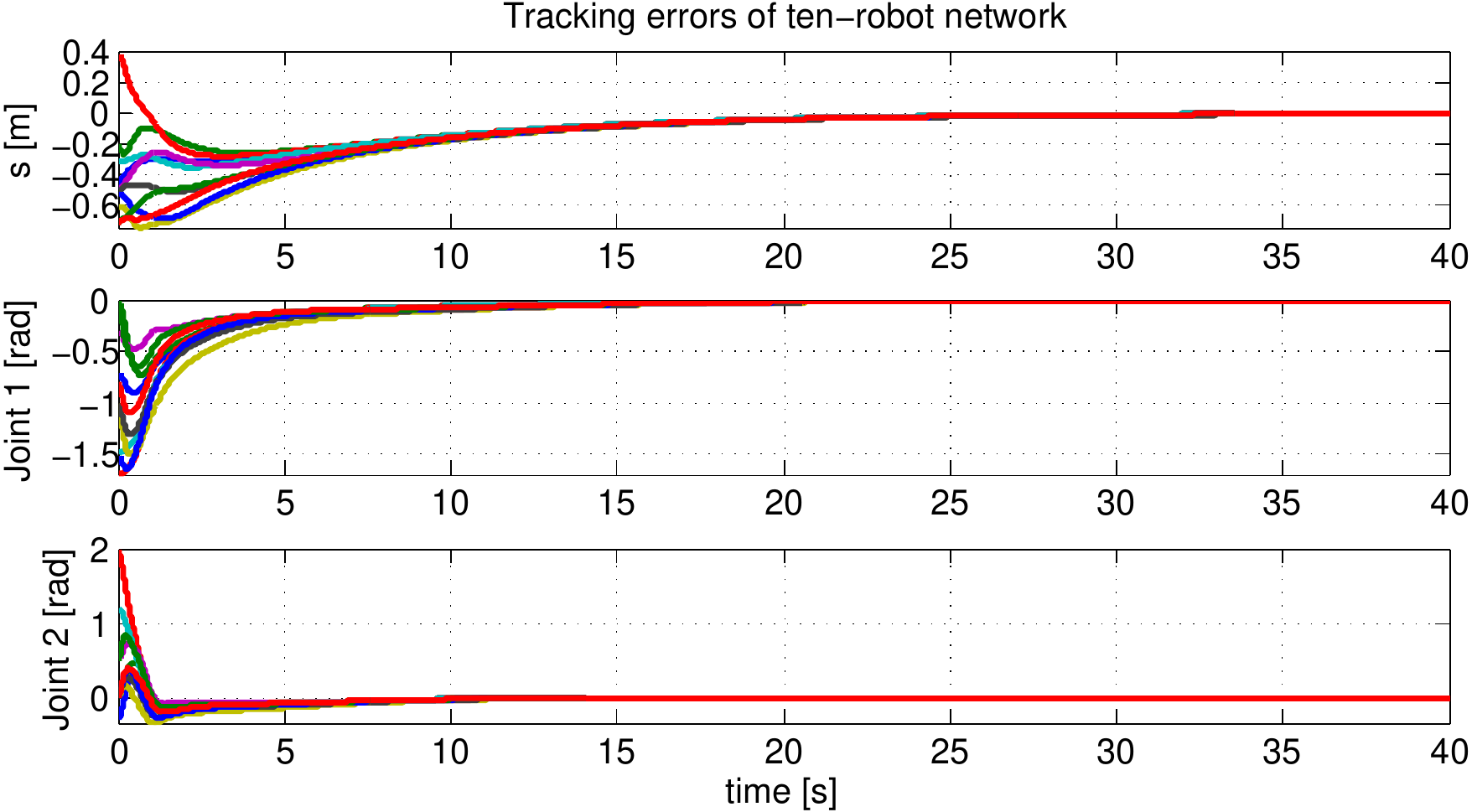}}
 \end{subfigmatrix}
 \caption{Synchronization of ten robots on three heterogeneous networks}\label{ten_robot_sim}
\end{figure}
\subsection{Simulation of Concurrent Synchronization for Ten Robots}
Let us consider three heterogeneous networks connected by feedback hierarchies shown in Fig.~\ref{ten_robot_sim}(a). If we view the reference trajectory inputs (dashed lines in the figure) as edges of the graph, as discussed in Section~\ref{sec:concurrent}, this graph is unbalanced. The first network (robot 1 through 4) is identical to the network of four 3-DOF robots in the previous section. All the physical parameters of the second network (robot 5 through 7) are twice larger than those of the first network while those of the third network are 1.5 times larger. Also, note that the third network is on an inline configuration by connecting the second feedback gain $\mathbf{K}_2$ of the control law in (\ref{tracking_controllerCH5}) back to itself. All the previous theorems still hold for such an inline configuration (see Corollary~\ref{CorollaryInline}). As seen in Fig.~\ref{ten_robot_sim}(b), the first, second, and third network individually synchronize robots within each network from arbitrary initial conditions. The first and second network also follow the reference trajectory command from the adjacent members. For the reference trajectory commands, we assume that we can send the position and velocity values ($\mathbf{q}_d,\mathbf{\dot q}_d)$. Then, $\mathbf{\ddot q}_d$ is estimated by a high pass filter. Eventually, all the robots exponentially synchronize and then follow the desired trajectory staying together. This concurrent synchronization can be used to construct a complex and large-scale dynamic network consisting of an arbitrary number of heterogeneous robots and networks.
\subsection{Simulation of Adaptive Synchronization}\label{sec:adaptive_sim}
\begin{figure}
 \begin{subfigmatrix}{1}
  \subfigure[with asymptotically stable tracking]{\includegraphics[width=3.3in]{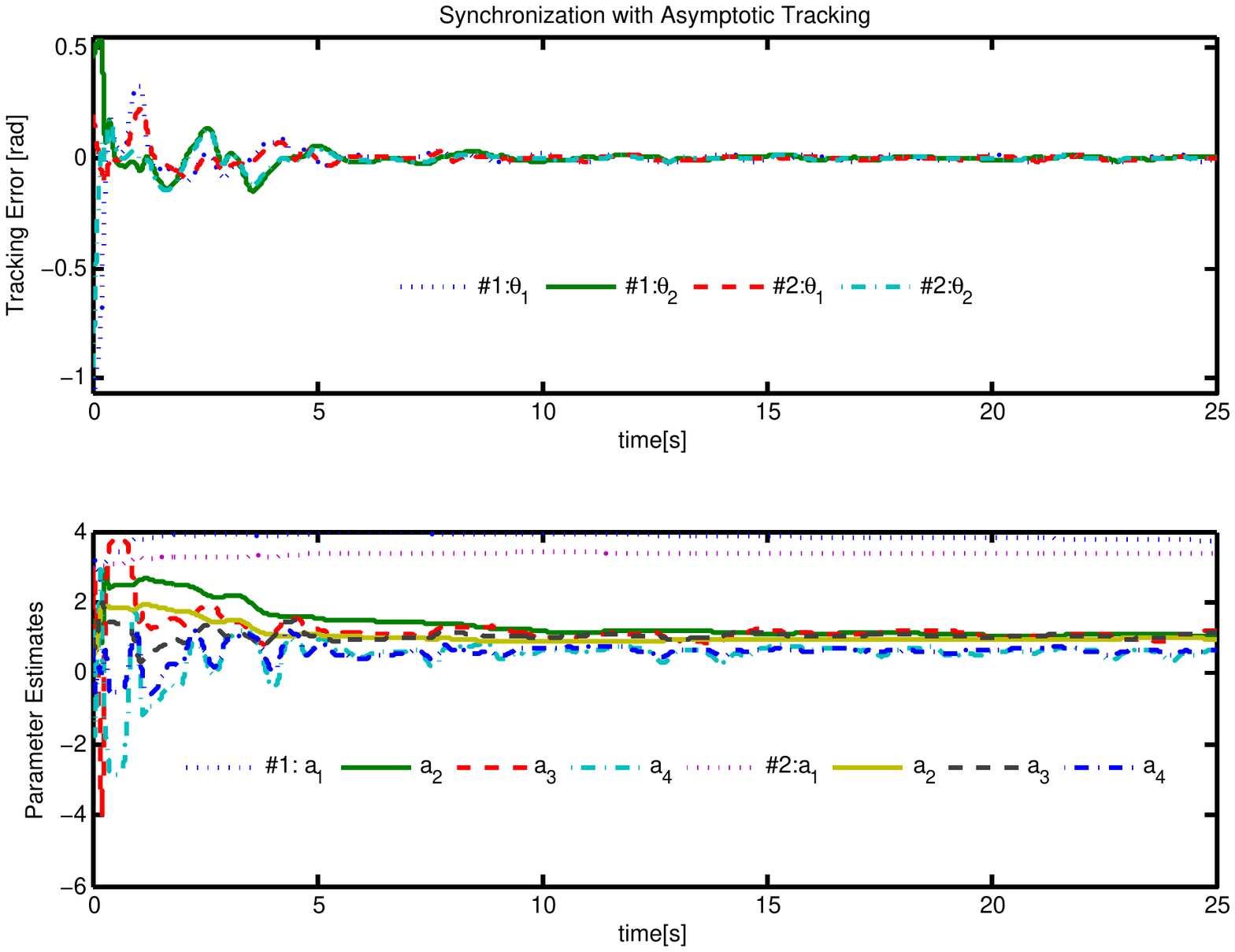}}
  \subfigure[with indifferent tracking]{\includegraphics[width=3.3in]{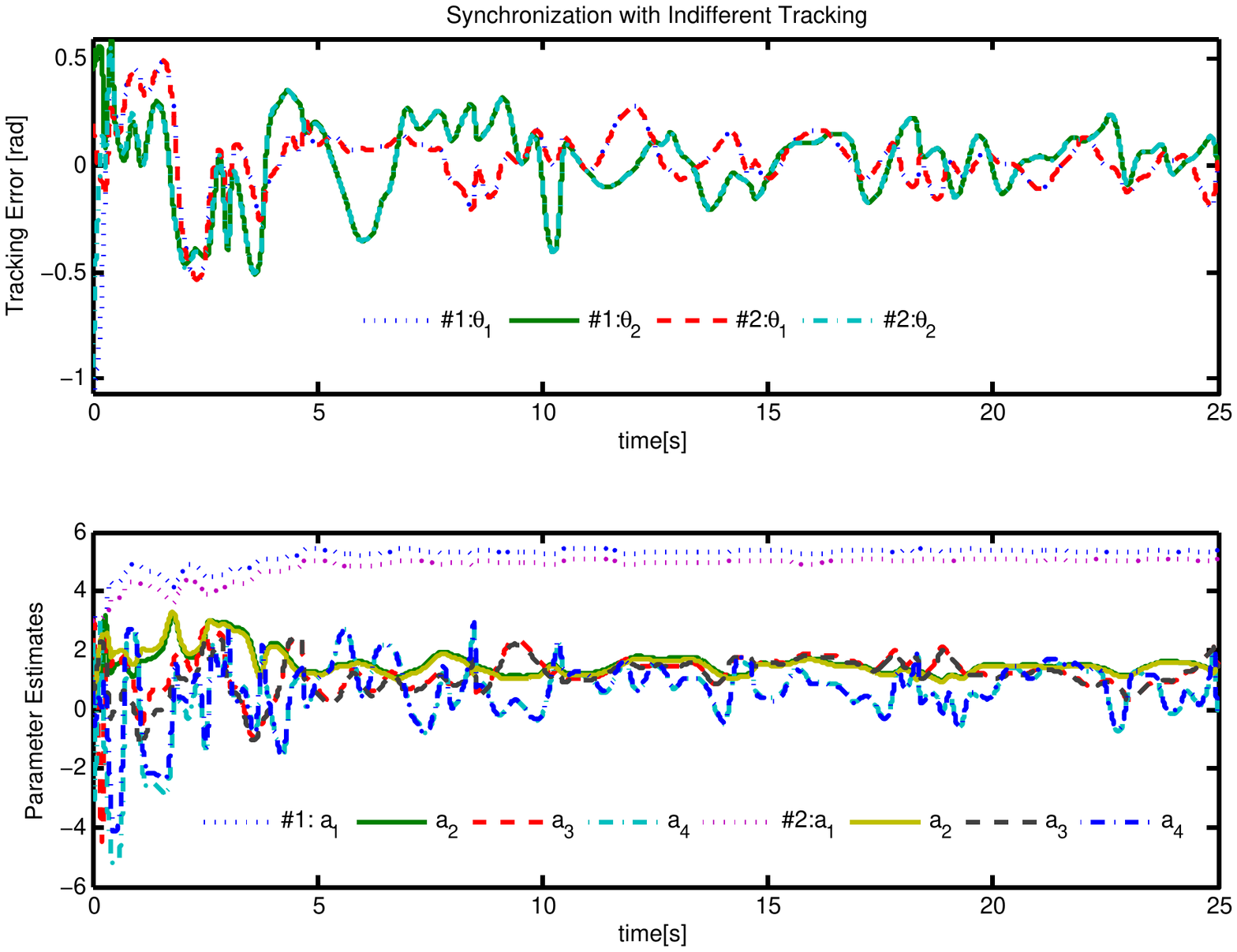}}
 \end{subfigmatrix}
 \caption{Adaptive synchronization of two two-link manipulator robots}\label{adaptive_sim}
\end{figure}

We assess the effectiveness of the proposed adaptive synchronization control law presented in Section~\ref{sec:adaptive_sync}. Consider a two-robot network, comprised of the two-link manipulator robots given in page 396 of \cite{Ref:Slotine}. The desired trajectory for the first joint is $\theta_{1d}(t)=\sin{\pi t}$ and for the second joint is $\theta_{2d}(t)=2(1-\cos{0.6\pi t})$. The initial parameter estimates for both robots are defined as $\mathbf{\hat{a}}(0)=(3,1,1,1)^T$, and $\mathbf{\Gamma}=\text{diag}(0.03,0.05,0.1,0.3)$. Figure~\ref{adaptive_sim}(a) shows the synchronization of the two robots with stable tracking by $\mathbf{K}_1=20\mathbf{I}$, $\mathbf{K}_2=15\mathbf{I}$, and $\mathbf{\Lambda}=10\mathbf{I}$. Hence, the synchronization gain $\mathbf{K}_1+\mathbf{K}_2$ is larger than the tracking gain $\mathbf{K}_1-\mathbf{K}_2$. Note that the synchronization of the tracking errors implies the synchronization of the state variables (i.e., joint 1 and joint 2). As discussed before, the synchronization occurs faster than the tracking. This simulation result indicates that the proposed adaptive control law can be used to synchronize motions of robots with unknown physical parameters.

In contrast, Figure~\ref{adaptive_sim}(b) shows the synchronization with indifferent tracking by the gains of $\mathbf{K}_1=\mathbf{K}_2=20\mathbf{I}$ and $\mathbf{\Lambda}=10\mathbf{I}$. The tracking dynamics then have zero gain (indifferent). So the two robots synchronize asymptotically, while the tracking errors remain within a finite ball. In both cases (Figs.~\ref{adaptive_sim}(a) and (b)), the proposed adaptive control law ensures neither the asymptotic convergence nor the synchronization of the physical parameter estimates $\mathbf{\hat{a}}_i$, unless the persistency of excitation condition is met. Nevertheless, the robots synchronize asymptotically from any initial condition in the presence of such parametric uncertainty. This adaptive version makes the proposed synchronization framework more practical. It should also be obvious to readers that the proposed control law can be straightforwardly extended to robust adaptive control schemes based on sliding control~\cite{Ref:Slotine}.

\section{Conclusions}
We have presented the new synchronization tracking control
law that can be directly applied to cooperative control of
multi-robot systems and oscillation synchronization in robotic
manipulation and locomotion. We have also shown that complex dynamic networks can be constructed by exploiting concurrent synchronization such that multiple groups of fully synchronized elements coexist. The proposed decentralized control law,
which requires only local coupling feedback for global exponential
convergence, eliminates both the all-to-all coupling and the
feedback of the acceleration terms, thereby reducing communication
burdens and complexity. Furthermore, in contrast with prior
work which used simple single or double integrator models, the
proposed method permits highly nonlinear systems. Providing exact nonlinear stability results constitutes one of the main contributions of this paper; global and exponential stability of the
closed-loop system has been derived by using contraction theory.
Contraction analysis, overcoming a local result of Lyapunov's
indirect method, yields global results based on differential
stability analysis. While we have focused on the mutual synchronization problem where synchronization and trajectory tracking take place simultaneously, the proposed method has also been shown to be a generalization of the average consensus problem that does not address trajectory tracking. It has been emphasized that there exist multiple time scales in the closed-loop systems: the faster convergence rates represents the transient boundary layer dynamics of synchronization while the slower rate determines how fast the synchronized systems track the common reference trajectory. Exponential synchronization with a faster convergence
rate enables reduction of multiple dynamics into a simpler form,
thereby simplifying the stability analysis. 

Simulation results show the effectiveness of the proposed control strategy. The proposed
bi-directional coupling has also been generalized to permit
partial-state coupling and uni-directional coupling. Further
extensions to proportional-derivative coupling, adaptive
synchronization, time-varying network topology, and concurrent synchronization of multiple heterogeneous networks on unbalanced graphs exemplify the benefit
of the proposed approach based on contraction
theory.
\section*{Acknowledgements}
This paper benefited from constructive comments from the associate editor, anonymous reviewers, Dr. Keehong Seo, and Prof. David W. Miller at MIT.

\begin{thebibliography}{99}
\bibitem{Ref:Arimoto1}S. Arimoto, F. Miyasaki, H. G. Lee, and S. Kawamura, ``Revival of Lyapunov's direct method in robot control and design,'' in \emph{Proc. American Control Conf.}, Atlanta, GA, June 1988, pp. 1764--1769.
\bibitem{Ref:Arimoto2}S. Arimoto, \emph{Control Theory of Nonlinear Mechanical Systems: A Passivity-based and circuit-theoretic approach}, Oxford University Press, Oxford, UK, 1996.
\bibitem{Ref:Spong_network}
N. Chopra and M. W. Spong, ``On synchronization of networked
passive systems with time delays and application to bilateral teleoperation,''
\emph{Proceedings of the SCIE Annual Conference}, Okayama, Japan, August 8-10, 2005, pp. 3424-–3429.
\bibitem{Ref:Spong_Sync}
N. Chopra and M. W. Spong, ``Passivity-based control of multi-agent systems,'' in \emph{Advances in Robot Control: From everyday physics to human-like movements}, Sadao Kawamura and Mikhail Svinin, Editors, Springer-Verlag, Berlin, 2006, pp. 107--134.
\bibitem{Ref:SpectralGraph}
F. Chung, \emph{Spectral Graph Theory}, Number 92 in CBMS Regional
Conference Series in Mathematics, American Mathematical Society,
1997.
\bibitem{Ref:sjchung_CDC}
S.-J. Chung and J.-J. E. Slotine, ``Cooperative robot control and synchronization of Lagrangian systems,'' \emph{Proc. 46th IEEE Conf. on Decision and Control}, New Orleans, LA, Dec. 2007, pp. 2504--2509.
\bibitem{Ref:sjchung_PhD}
S.-J. Chung, ``Nonlinear control and synchronization of
multiple Lagrangian systems with application to tethered formation
flight spacecraft,'' Doctor of Science Thesis, Dept. of Aeronautics and
Astronautics, Massachusetts Inst. of
Technology, Cambridge MA, 2007.
\bibitem{Ref:tether}
S.-J. Chung, J.-J. E. Slotine, and D. W. Miller, ``Nonlinear model
reduction and decentralized control of tethered formation flight,''
\emph{Journal of Guidance, Control, and Dynamics}, Vol. 30, No. 2, pp. 390--400, March--April 2007.
\bibitem{Ref:FormationFlying_CHUNG}
S.-J. Chung, U. Ahsun, and J.-J. E. Slotine, ``Application of synchronization to formation flying spacecraft: Lagrangian approach,'' accepted for publication, \emph{Journal of Guidance, Control and Dynamics}, to be published.
\bibitem{Demi}
B. P. Demidovich, ``Dissipativity of a nonlinear system of
differential equations,'' {\it ser. matematika mehanika, part I {\bf
N.6}, pp. 19--27; part II {\bf N.1}, pp.3--8} Vestnik Moscow State
University, 1962.
\bibitem{Ref:robot_formation}
R. Fierro, P. Song, A. Das, and V. Kumar, ``Cooperative control
of robot formations,'' in \emph{Cooperative Control and Optimization:
Series on Applied Optimization}, Kluwer Academic Press, pp. 79--93, 2002.

\bibitem{Ref:Gerard_Slotine}
L. Gerard and J.-J. E. Slotine, ``Neuronal networks and controlled symmetries, a generic framework,'' arXiv:q-bio/0612049v2 [q-bio.NC].
\bibitem{Hart}
P. Hartmann, \emph{Ordinary Differential Equations} {\em John Wiley
\& Sons, New York}, 1964.
\bibitem{IhleArcak}
I.-A. F. Ihle, M. Arcak, and T. I. fossen, ``Passivity-based designs for synchronized path-following,'' \emph{Automatica}, Vol. 43, No. 9, pp. 1508--1518, 2007.
\bibitem{Jadbabaie}
A. Jadbabaie, J. Lin, and A. S. Morse, ``Coordination of groups
of mobile autonomous agents using nearest neighbor rules,'' \emph{IEEE
Trans. on Automatic Control}, Vol. 48, No. 6, pp. 988--1001, June 2003.
\bibitem{Jouffroy}
J. Jouffroy and J.-J. E. Slotine, ``Methodological remarks on contraction theory,'' \emph{IEEE Conf. on Decision and Control}, Atlantis, Paradise Island, Bahamas, 2004.
\bibitem{Khalil:2002}
H. K. Khalil, \emph{Nonlinear Systems}, 3rd Ed.,Prentice Hall, Upper
Saddle River, NJ, 2002, pp. 339--350.
\bibitem{Krstic}
K. Krstic, I. Kanellakopoulos, and P. Kokotovic, \emph{Nonlinear and Adaptive control Design},
John Wiley and Sons, New York, NY, 1995, pp. 21--29.
\bibitem{Ref:passive_decomp1}
D. Lee and P. Y. Li, ``Formation and maneuver control of multiple
spacecraft,'' \emph{Proc. of the 2003 American Control Conf.},
Vol. 1, June 2003, pp. 278--283.
%
\bibitem{Ref:passive_decomp2}
D. Lee and M. W. Spong, ``Stable flocking of multiple inertial
agents on balanced graph,'' \emph{Proc. of the 2006 American Control
Conf.}, Vol. 52, June 2006, pp. 1469--1475.

\bibitem{Ref:Leonard_leader}
N. E. Leonard and E. Fiorelli, ``Virtual leaders, artificial
potentials and coordinated control of groups,'' \emph{Proc. 40th IEEE Conf. on Decision and Control}, 2001, pp. 2968--2973.
\bibitem{Lew}
D. C. Lewis, ``Metric
properties of differential equations,'' {\it American Journal of
Mathematics,} Vol. 71, pp. 294--312, 1949.

\bibitem{Ref:mobile}
Z. Lin, M. Broucke, and B. Francis, ``Local control strategies
for groups of mobile autonomous agents,'' \emph{IEEE Trans. on Automatic
Control}, Vol. 49, No. 4, pp. 622--629, Apr. 2004.

\bibitem{Ref:contraction1}
W. Lohmiller and J.-J. E. Slotine, ``On contraction analysis for
nonlinear systems,'' \emph{Automatica}, Vol. 34, No. 6, pp. 683–-696, 1998.

\bibitem{Ref:Slotine_highorder}
W. Lohmiller and J.-J. E. Slotine, ``High-order nonlinear
contraction analysis,'' MIT Nonlinear Systems Laboratory (NSL) Report, NSL-050901.
\bibitem{Ref:Mesbahi}
M. Mesbahi and F. Y. Hadaegh, ''Formation flying of multiple
spacecraft via graphs, matrix inequalities, and switching,'' \emph{Journal of Guidance, Control, and Dynamics}, Vol. 24, No. 2, pp. 369--377, 2001.
\bibitem{Ref:tele}
G. Niemeyer and J.-J. E. Slotine, ``Telemanipulation with time
delays,'' \emph{Int. J. of Robotics Research}, Vol. 23, No. 9, pp. 873--890, 2004.
\bibitem{Ref:Murray}
R. Olfati-Saber and R. M. Murray, ``Consensus problems in networks
of agents with switching topology and time-delays,'' \emph{IEEE Trans. on Automatic Control}, Vol. 49, No. 9, pp. 1520--1533, Sep. 2004.
\bibitem{CLF:2002}
P. \"Ogren, M. Egerstedt, and X. Hu, ``A control Lyapunov function
approach to multiagent coordination,'' \emph{IEEE Trans. on
Robotics and Automation}, Vol. 18, No. 5, pp. 847--851, Oct. 2002.
\bibitem{Ref:multi_coord}
P. \"Ogren, E. Fiorelli, and N. E. Leonard, ``Cooperative control
of mobile sensor networks: adaptive gradient climbing in a
distributed environment,'' \emph{IEEE Trans. on Automatic Control},
Vol. 49, No. 8, pp. 1292--1302, 2004.
\bibitem{PaleyLeonard}
D. A. Paley, N. E. Leonard, R. Sepulchre, D. Grunbaum, and J. K. Parrish, ``Oscillator models and collective motion: spatial patterns in the dynamics of engineered and biological networks,'' \emph{IEEE Control Systems Magazine}, Vol. 27, No. 4, pp. 89--105, August, 2007

\bibitem{Ref:contraction_sync}
Q.-C. Pham and J.-J. E. Slotine, ``Stable concurrent synchronization
in dynamic system networks,'' \emph{Neural Networks}, Vol. 20, No. 1, pp. 62--77, 2007.
\bibitem{Ref:sync_Pitti}
A. Pitti, M. Lungarella, and Y. Kuniyoshi, ``Synchronization:
adaptive mechanism linking internal and external dynamics,'' \emph{Proc.
of the 6th Int. Workshop on Epigenetic Robotics}, Paris, France, Sept. 2006.
\bibitem{Ref:WeiRen1}
W. Ren, ``Multi-vehicle consensus with a time-varying reference state,'' \emph{Systems \& Control Letters}, Vol. 56, No. 7-8, pp. 474--483, July, 2007.
\bibitem{Ref:WeiRen2}
W. Ren, R. W. Beard, and E. Atkins, ``Information consensus in multivehicle cooperative control: collective group behavior through local interaction,'' \emph{IEEE Control Systems Magazine}, Vol. 27, No. 2, pp. 71--82, April, 2007.
\bibitem{Ref:robot_sync}
A. Rodriguez-Angeles and H. Nijmeijer, ``Mutual synchronization of robots via
estimated sate feedback: a cooperative approach,'' \emph{IEEE Trans.
on Control Systems Technology}, Vol. 12, No. 4, pp. 542--554, July 2004.
\bibitem{SeoA}K. Seo, S.-J. Chung, and J.-J. E. Slotine, ``CPG-based control of a turtle-like underwater vehicle,'' \emph{Proc. of the Robotics: Science and Systems (RSS)}, Switzerland, June 2008.
\bibitem{Ref:Slotine}
J.-J. E. Slotine and W. Li, \textit{Applied Nonlinear Control}, Prentice Hall,. New
Jersey, 1991.
\bibitem{Ref:contraction2}
J.-J. E. Slotine, ``Modular stability tools for distributed
computation and control,'' \emph{Int. J. Adaptive Control and Signal
Processing}, Vol. 17, No. 6, pp. 397--416, 2003.
\bibitem{Ref:contraction5}
J.-J. E. Slotine and W. Lohmiller, ``Modularity, evolution, and the binding problem: a view from stability theory,'' \emph{Neural Networks}, Vol. 14, No. 2, pp. 137--145, 2001.
\bibitem{Ref:contraction_robot}
J.-J. E. Slotine, W. Wang, and K. El Rifai, ``Synchronization in
networks of nonlinearly coupled continuous and hybrid oscillators,'' \emph{Proc.
of the 6th International Symposium on Mathematical Theory of Networks and Systems (MTNS 2004)}, Katholieke Universiteit Leuven, Leuven, Netherlands, Jul. 5--9, 2004.
\bibitem{Ref:Strang}
G. Strang, \emph{Introduction to Applied Mathematics},
Wellesley-Cambridge Press, Wellesey, MA, 1986.
\bibitem{Sun1}
D. Sun and J. K. Mills, ``Adaptive synchronized control for coordination of multi-robot assembly tasks,'' \emph{IEEE Trans. on Robotics and Automation}, Vol. 18, No. 4, pp 498--510, Aug. 2002.
\bibitem{Sun2}
D. Sun, L. Ren, J. K. Mills, and C. Wang, ``A synchronous tracking control of parallel manipulators using cross-coupling approach,'' \emph{Int. J. of Robotics Research}, Vol. 25, No. 11, pp. 1137--1147, Nov. 2006.
\bibitem{Ref:contraction3}
W. Wang and J.-J. E. Slotine, ``On partial contraction analysis for
coupled nonlinear oscillators,'' \emph{Biological Cybernetics}, Vol. 92, No. 1, pp. 38–-53, 2005.
\bibitem{Ref:contraction_timedelay}
W. Wang and J.-J. E. Slotine, ``Contraction analysis of time-delayed communications and group cooperation,'' \emph{IEEE Trans. on Automatic Control}, Vol. 51, No. 4, pp. 712--717, Apr. 2006.
\bibitem{Ref:contraction_leader}
W. Wang and J.-J. E. Slotine, ``A theoretical study of different
leader roles in networks,'' \emph{IEEE Trans. on Automatic Control}, Vol. 51, No. 7, pp. 1156--1161, 2006.
\end{thebibliography}
\end{document}